\documentclass{amsart}
\usepackage[latin1]{inputenc}
\usepackage{amssymb}
\usepackage{amsmath}
\usepackage{enumerate}
\usepackage{epsfig}
\usepackage{subfigure}
\usepackage[all]{xy}
\usepackage{longtable}

\newtheorem{theorem}{Theorem}[section]
\newtheorem{lemma}[theorem]{Lemma}
\newtheorem{proposition}[theorem]{Proposition}
\newtheorem{corollary}[theorem]{Corollary}
\newtheorem{condition}[theorem]{Condition}

\theoremstyle{definition}
\newtheorem{definition}[theorem]{Definition}
\newtheorem{example}[theorem]{Example}

\newtheorem{question}[theorem]{Question}

\theoremstyle{remark}
\newtheorem{rem}[theorem]{Remark}
\theoremstyle{remark}

\numberwithin{equation}{section}

\newcommand{\rk}{\mbox{rank}}

\newcommand{\ra}{\rightarrow}

\newcommand{\Z}{\mathbb{Z}}
\newcommand{\Q}{\mathbb{Q}}

\newcommand{\N}{\mathbb{N}}

\newcommand{\Aut}{{\operatorname{Aut}}}
\newcommand{\Nef}{{\operatorname{Nef}}}
\newcommand{\Pic}{\operatorname{Pic}}

\makeatletter
\def\blfootnote{\xdef\@thefnmark{}\@footnotetext}

\begin{document}

\subjclass[2010]{Primary 14J28; Secondary  14E30.}
\keywords{K3 surfaces, Mori Dream Spaces, Extremal Contractions.\\
The author is partially supported by PRIN 2010--2011 ``Geometria delle variet\`a algebriche" and FIRB 2012 ``Moduli Spaces and their Applications".}

\title[Mori Dream Spaces extremal contractions of K3 surfaces]{Mori Dream Spaces extremal contractions of K3 surfaces}
\author{Alice Garbagnati}
\address{Alice Garbagnati, Dipartimento di Matematica, Universit\`a di Milano,
  via Saldini 50, I-20133 Milano, Italia}
\email{alice.garbagnati@unimi.it}

\begin{abstract}
We will give a criterion to assure that an extremal contraction of
a K3 surface which is not a Mori Dream Space produces a singular
surface which is a Mori Dream Spaces. We list the possible
N\'eron--Severi groups of K3 surfaces with this property and an extra geometric condition such that the Picard number is greater then or equal to 10. We give a
detailed description of two geometric examples for which the
Picard number of the K3 surface is 3, i.e. the minimal possible in
order to have the required property. Moreover we observe that
there are infinitely many examples of K3 surfaces with the
required property and Picard number equal to 3.
\end{abstract}

\maketitle
\section{Introduction}

The Mori Dream Spaces are projective varieties for which the
Minimal Model Program can be applied successfully, i.e. the
necessary flops and contractions exist and the program terminates.
They were introduced in \cite{HK} where it is also proved that the
property of being a Mori Dream Space is essentially equivalent to
the finite generation of the Cox ring of the variety.

In this paper we discuss the following situation: let us consider
a variety which is not a Mori Dream Space and let us contract
exactly one extremal ray of its effective cone. Is it possible to
obtain in such a way a Mori Dream Space? The answer to this
question is known to be yes, so the property "to be a Mori Dream
Space" is yet not preserved by resolution of simple singularities, even if the resolution is a crepant resolution.
The examples which allow one to give a positive answer to this
question are provided in \cite[Theorem 5.2]{O}. These examples are
constructed as follows: let $S$ be a K3 surface which is a Mori
Dream Space. Then the singular variety $V:=(S\times
S)/\mathfrak{S}_2$, where $\mathfrak{S}_2$ is the symmetric group,
is known to be a Mori Dream Space. The variety $V$ is singular
along the diagonal and there exists an extremal crepant resolution
of $V$ which is the Hyperk\"ahler variety $H:=Hilb^2(S)$. The
variety $V$ is obtained from $H$ by contracting the exceptional
divisor over the diagonal, so $V$ is an extremal contraction of
$H$. In \cite{O} Oguiso gives examples of K3 surfaces $S$ such
that $S$ and $V$ are Mori Dream Spaces but $H$ is not. The proof
of the fact that $S$ and $V$ are a Mori Dream Space and $H$ is not
is based on the computation of the automorphism groups of these
varieties. Indeed, it is possible to construct K3 surfaces $S$
such that $S$ and $V$ have a finite automorphism groups, but the
automorphism group of $H$ is not finite.

The relation among the finiteness of the automorphism group and
the finite generation of the Cox ring (so the property of being a
Mori Dream Space), is deep as shown for example by the following
result, by Artebani, Hausen and Laface, \cite{AHL}, which will be
be fundamental in this paper: a K3 surface is a Mori Dream Space
if and only if its automorphism group is finite.

In view of the result by Oguiso on the varieties $H$ and $V$, it
is quite natural to ask if there exist similar results in lower
dimension, in particular for K3 surfaces. So the aim of this note
is to positively answer to the following question:

\begin{question}\label{question} Are there K3 surfaces $X$ which are not Mori Dream Spaces and such that the surface $X'$ obtained by contracting exactly one extremal ray of $X$ is a Mori Dream Space?
\end{question}

We will answer to the question \ref{question} both by providing an
abstract algorithm which produces the N\'eron Severi groups of
admissible examples and describing the geometry of some of them.

Since $X$ is a K3 surface, it is by definition smooth and minimal
and so the extremal contraction of $X$ produces a singular surface
$X'$. From this point of view, $X$ is a crepant extremal
resolution of the singular surface $X'$, and we are constructing a
crepant extremal resolution of a Mori Dream Space which is not a
Mori Dream Space. The same is true for the Oguiso's examples,
indeed the variety $H$ is crepant extremal resolution on $V$.

In Section \ref{sec: criterion} we prove that under a geometrical condition on the smooth rational curves on $X$, the singular surface
$X'$ is a Mori Dream Space if and only if the K3 surface $Y$ whose
N\'eron--Severi group is isometric to the one of $X'$ is a Mori
Dream Space.  In Section \ref{sec: pairs X and X'} we give a
criterion to find  the N\'eron--Severi group of the K3 surfaces
$X$ such that there exists an extremal contraction producing a
singular surface $X'$ which is a Mori Dream Space. Here we use the
results of \cite{AHL} on K3 surfaces which are Mori Dream Spaces
and the results of \cite{K} on K3 surfaces with finite
automorphism group. The main result of the Section \ref{sec: pairs
X and X'} is Theorem \ref{theo: the pairs X,X'} where we prove
that, under a geometric condition on the rational curves on $X$, if $X$ is a K3 surface which is not a Mori Dream Space, but
which admits an extremal contraction producing a singular surface
which is a Mori Dream Space, then $\rho(X)\geq 3$ and $\rho(X)\neq
19$. Moreover, we classify all the admissible N\'eron--Severi
groups of $X$ if $\rho(X)\geq 10$ (the reason of this bound is
only computational, and is essentially due to the fact that in
\cite{K} the same bound is considered). Section \ref{sec:
geometric example} is the geometric heart of this paper. Here we
give two examples of K3 surfaces $X$ which are not Mori Dream
Spaces, but which admit an extremal contraction producing a
singular surface which is a Mori Dream Space. In both the cases
$\rho(X)=3$ and the N\'eron--Severi groups of the surfaces
obtained after the contraction are the same. We will show that the
automorphism group of the K3 surface $X$ is infinite and we show
that some automorphisms do not descend to automorphisms of the
singular model. In Section \ref{sec: other examples} we provide
other examples, giving some geometric details. In particular we
will show that there exists an infinite number of K3 surfaces of
Picard number 3 which are not Mori Dream Spaces but such that the
singular surface obtained by one extremal contraction is a Mori
Dream Space, see Proposition \ref{prop: infinite pairs with K3
with rank 3 not Mori} for a more precise statement.

{\bf Acknowledgements}: I would like to thank Professor Oguiso for
suggesting me this problem, for our useful discussions and for the
kind hospitality when I visited Osaka's University. I would also
like to thank Professor van Geemen for several suggestions. 
I would like to show my gratitude to the anonymous referee who found and corrected a mistake in the preliminary version and for essential comments to this paper.

\section{A criterion to conclude that the contraction of an extremal ray of a K3 surface produces a Mori Dream Space}\label{sec: criterion}

\subsection{The K3 surfaces $X$ and $Y$ and the singular surface $X'$}
\begin{definition} A K3 surface is a complex compact surface with
trivial canonical bundle and trivial irregularity.

We denote by $\Lambda_{K3}$ the unique even unimodular lattice
with signature $(3,19)$ and we observe that if $S$ is a K3
surface, then $H^2(S,\Z)\simeq \Lambda_{K3}$.

In the following for each lattice $L$ and each element $l\in L$ we
denote by $l^{\perp_L}$ the sublattice of $L$ whose elements are
orthogonal to $l$, i.e. $l^{\perp_L}:=\{k\in L\mbox{ such that
}kl=0\}$.\end{definition}

Let $X$ be an algebraic K3 surface (i.e. a K3 surface with an ample line bundle). The N\'eron--Severi group of $X$ is a sublattice of $\Lambda_{K3}$ and it is hyperbolic (i.e. its signature is $(1,\rho(X)-1)$). We assume that $\rho(X):=\rk(NS(X))\geq 2$. We
denote by $L$ the abstract lattice such that $L\simeq NS(X)$. Let
us assume that there exists a smooth rational curve $N$ on $X$. We
denote by $N$ also the class of the curve $N$ in $NS(X)\simeq L$
and we recall that $N^2=-2$. We denote by $M$ the sublattice of
$NS(X)\simeq L$ which is orthogonal to $N$, i.e. $M:=N^{\perp_L}$.
Since $L$ is a hyperbolic lattice and $\langle N \rangle$ is a
negative definite lattice, the lattice $M$ is a hyperbolic
lattice. Moreover, the lattice $M$ is primitively embedded in $L$,
by definition. Since there exists a primitive embedding of $L$ in
$\Lambda_{K3}$, there exists a primitive embedding of $M$ in
$\Lambda_{K3}$. As a consequence of the surjectivity of the period
map for K3 surfaces, there exists at least a K3 surface (and in
fact infinitely many) whose N\'eron--Severi group is isometric to
$M$. Let $Y$ be a K3 surface such that $NS(Y)\simeq M$. We observe
that $\rho(Y)=\rho(X)-1$.

The curve $N$ is a smooth irreducible curves on $X$, so it determines a wall of the chamber of the positive cone which coincides with the ample cone of $X$. Since the set of $\Q$-divisors is dense in the positive cone, there exists at least one pseudoample divisor which is on the wall determined by $N$, i.e. there exists a pseudoample divisor orthogonal to $N$. Since the set of wall is locally finite, there exists at least one pseudoample divisor which is on the wall determined by $N$, but not on other walls, i.e. there exists a pseudoample divisor, say $H_X$, which is orthogonal to $N$, but not to any other smooth rational curve on $X$. Since the divisor $H_X$ is orthogonal to $N$, it is represented by a vector in $M\subset L\simeq NS(X)$. The same vector represent a divisor $H_Y\in M\simeq NS(Y)$ on $Y$, which is a divisor with a positive self intersection and with a non zero-intersection with all the irreducible classes on $M$ with self intersection $-2$. Hence we can chose the ample cone of $Y$ to be the chamber of the positive cone of $Y$ which contains $H_Y$. 
To recap we have a K3 surface $X$ on which we fixed a smooth rational curve $N$ and a pseudoample divisor $H_X$ (which is orthogonal to $N$ and not to other irreducible smooth rational curves on $X$). We associate to $(X,N,H_X)$ a K3 surface $Y$ such that $NS(Y)\simeq M\simeq N^{\perp_{NS(X)}}$ and $H_Y$ (the divisor represented in $M$ by the same vector which represents $H_X$) is ample.

Let us denote by $\mathcal{B}_M=\{B_1,\ldots B_{\rk(M)}\}$ a
$\Z$-basis of the abstract lattice $M$. The set
$\mathcal{B}_L=\{B_1,\ldots B_{\rk(M)}, N\}$ is a $\Q$-basis of
the lattice $L$, i.e. every element in $L$ is a linear combination
with rational coefficients of $\{B_1,\ldots B_{\rk(M)}, N\}$. For
certain specific lattice $L$ the set $\mathcal{B}_L$ is also a
$\Z$-basis for $L$, but this can not be guarantee in a general
context. However, every element in $L$ which has a trivial
intersection with $N$ is a linear combination with integer
coefficients of $\{B_1,\ldots B_{\rk(M)}\}$.

Let $\phi:X\rightarrow X'$ be the map which contracts the curve
$N$ to a point. Then $X'$ is a singular variety with exactly one
singular point, of type $A_1$, which is the point $\phi(N)$.

The surface $X'$ is a normal surface and its N\'eron--Severi group
is isometric to $M$. So, the N\'eron--Severi group  $NS(X')$ and
the N\'eron--Severi group $NS(Y)$ are isometric, and both
isometric to $M$.

By construction, we fixed a specific primitive embedding of $M$ in
$NS(X)$, which is in fact a marking. Let us denote this embedding
by $f_X:M\ra NS(X)$. So $NS(X')$ and $NS(Y)$
are identified with $N^{\perp_L}$, and the marking $f_X$ induces
the isometries $f_{X'}:M\ra NS(X')$ by $f_Y:M\ra NS(Y)$. For each
vector $v\in M$ we denote by $D_X:=f_X(v)$, $D_{X'}:=f_{X'}(v)$
and $D_{Y}:=f_{Y}(v)$. We say that the divisors $D_X$, $D_{X'}$
and $D_Y$ are associated if they are the images (for the map
$f_X$, $f_{X'}$ and $f_Y$ respectively) of the same vector $v\in
M$. Let $D_X$ be a divisor on $X$ such that $D_X=f_{X}(v)$ for a
certain $v\in M$. Then $D_XN=0$ and so the associated vectors are
$D_{X'}=f_{X'}(f_X^{-1}(D_X))$ and $D_{Y}=f_{Y}(f_X^{-1}(D_X))$.

\subsection{The Nef cones of $Y$ and $X'$}
\begin{lemma}\label{lemma: from X to X'} Let $D_X$ be a divisor on $X$ such that $D_X$ is nef and $D_XN=0$. The divisor $D_{X'}\in NS(X')$ associated to $D_X$ is nef.\end{lemma}
\proof The Lemma follows directly by \cite[Example 1.4.4,ii)]{L}, but here we give a direct proof.
We denote by $v_D$ the vector in $M$ such that $f_X(v_D)=D_X$, and we recall that $D_{X'}=f_{X'}(v_D)$.
Since $D_X$ is nef, for every curve $C_X\subset X$, $D_XC_X\geq
0$. Let us now consider a curve $C_{X'}\subset X'$.
The strict transform of $C_{X'}$ on $X$
is a curve, $C_{X}$, whose class is $\alpha v_C+\eta N$, where $v_C\in f_X(M)$ and $\alpha,\eta\in \Q$. Hence the class of $C_{X'}$ is the $\Q$-divisor $\alpha f_{X'}(v_C)$. So we
have
$$D_{X'}C_{X'}=v_D\left(\alpha v_C\right)=v_D\left(\alpha v_C+\eta N\right)=D_XC_X\geq 0,$$
where we used that $f_{X}$ and $f_{X'}$ are isometries and the fact that $v_DN=0$.
Hence $D_{X'}$ has
a non negative intersection with all the curves in $X'$. We
conclude that $D_{X'}$ is nef.
\endproof

\begin{lemma} \label{lemma: from X to Y} Let $D_X$ be a divisor on $X$ such that $D_X$ is nef and $D_XN=0$.
The divisor $D_Y\in NS(Y)$ associated to $D_X$ is nef.\end{lemma}
\proof We consider a divisor $D_X$ which is nef. Let us consider the associated divisor $D_{Y}\in
NS(Y)$. It suffices to show that for every effective $(-2)$-class
$R_Y\in NS(Y)$, $R_YD_Y\geq 0$. Let $R_Y$ be an effective
$(-2)$-class in $NS(Y)$. It is associated to a divisor $R_X\in
NS(X)$. Clearly $R_X^2=R_Y^2=-2$. Since $X$ is a K3 surface, by
the Riemann-Roch theorem there are only the following two
possibilities: either $R_X$ is effective or $-R_X$ is effective.
If $R_X$ is effective, then $R_XD_X\geq 0$. By definition,
$R_XN=0$ and $D_XN=0$ so $R_XD_X=R_YD_Y$. Hence, $R_YD_Y\geq 0$
and $D_Y$ is nef. So it now suffices to exclude that $-R_X$ is
effective. Let us consider the intersection $R_XH_X$, since both
$R_X$ and $H_X$ are contained in $M\subset NS(X)$, we have
$R_XH_X=R_YH_Y$, which is non negative, since $H_Y$ is ample and
$R_Y$ is effective. So $R_XH_X>0$ and $H_X$ is pseudoample, this
implies that $-R_X$ can not be effective.\endproof

\begin{lemma}\label{lemma:from X' to X} Let $D_{X'}\in NS(X')$ be a nef divisor.
Let $D_X$ be the divisor associated to $D_{X'}$. Then $D_X\in
NS(X)$ is nef (and $D_XN=0$).\end{lemma} \proof The Lemma follows directly by \cite[Example 1.4.4,i)]{L}, but here we give a direct proof.

In order to show that $D_X$ is nef we show that for every curve $C_X\in NS(X)$,
$D_XC_X\geq 0$. If $C_X=N$, then $D_XC_X=0$. We now assume that
$C_X\neq N$. So $\phi(C_X)\subset X'$ is a curve in $X'$ and $C_X$
is the strict transform of $\phi(C_X)$ with respect to the blow up
$\phi$. The class of $C_X$ is represented in $NS(X)$ by
the class $\alpha f_X(v_C)+\eta N$ for a certain $v_C\in M$ and $\alpha, \eta\in\Q$. Then $\phi(C_{X})$ is represented by $\alpha f_{X'}(v_C)$. Let us denote by $v_D\in M$ the
vector such that $D_X=f_X(v_D)$. Thus,
$$D_XC_X=f_X(v_D)(\alpha f_X(v_C)+\eta N)=\alpha v_Dv_C=f_{X'}(v_D)\left(\alpha f_{X'}(v_C)\right)=D_{X'}\phi(C_X)$$ where we used
that $f_X(v_D)N=0$ and that $f_X$ and $f_{X'}$ are isometries.
Since $D_{X'}$ is nef, $D_{X'}\phi(C_X)\geq 0$, so $D_XC_X\geq 0$
for every curve $C_X$ in $X$.\endproof

\begin{lemma}\label{lemma: from Y to X} Let $D_Y\in NS(Y)$ be a nef divisor on $Y$. Let $D_X\in NS(X)$ be the associated divisor on $X$. If there exists no a curve $B_X\subset X$ with self intersection $-2$ such that $B_XN=1$, then $D_{X}$ is nef (and $D_XN=0$).
\end{lemma}

\proof
It suffices to prove that $D_XC_X\geq 0$ for every $C_X$ curve on $X$ with self intersection $-2$. If $C_X=N$, then $D_XC_X=D_XN=0$. So we assume that $C_X\neq N$ and thus $C_XN\geq 0$. Let us denote by $c:=C_XN$. By the hypothesis we know that $c\neq 1$, hence either $c=0$, or $c\geq 2$.
The divisor $C_X$ can be written in $NS(X)\otimes \Q$ as $C_X=A_X-\frac{c}{2}N$, where $A_X\in NS(X)\otimes \Q$ is a $\Q$-divisor orthogonal to $N$. Since $D_XN=0$, $D_XC_X=D_XA_X$.
If $c\geq 2$, then $A_X^2=\left(C_X+\frac{c}{2}N\right)^2=-2-\frac{c^2}{2}+c^2\geq 0$. So any positive multiple of $A_X$ has a positive self intersection, in particular there exists a multiple $mA_X$, $m\in \N$ of $A_X$ which is a divisor in $NS(X)$. By Riemann--Roch theorem either $mA_X$ or $-mA_X$ is effective. Since $C_X$ is a curve, $C_XH_X> 0$ for the pseudoample divisor $H_X$, and thus $mA_XH_X>0$, so $mA_X$ is an effective divisor in $X$. Thus the $\Q$-divisor $A_Y$($=f_Y(f_X^{-1}(A_X)$) is an effective $\Q$-divisor on $Y$ (since $A_Y^2\geq 0$ and $A_YH_Y>0$). The divisor $D_Y$ is nef on $Y$, so $D_YA_Y\geq 0$. By $D_XC_X=D_XA_X=D_YA_Y\geq 0$ we conclude the proof in case $c\geq 2$.

If $c=0$, i.e. $CX_N=0$, $C_X$ is by definition an effective divisor on $X$ and the associated divisor $C_Y$ is an effective divisor on $Y$. So $D_YC_Y\geq 0$ and we conclude by $D_XC_X=D_YC_Y\geq 0$. 
\endproof

\begin{rem}\label{rem: the ccondition on Bx is necessarly}{\rm In Lemma \ref{lemma: from Y to X}, the assumption that there exists no a $(-2)$-curve $B_X$ such that $B_XN=1$ is essential. Indeed if the curve $B_X$ exists, then the divisor $D_X$ in the statement can not be nef. Because of the duality between the effective cone and the nef cone of a surface, to prove that $D_X$ can not be nef, it suffices
to produce an effective class $E_X$ on $X$ such that $E_XN=0$, but $E_Y$($=f_Y^{-1}(f_X(E_X))$) is not effective on $Y$. 

So we consider the curve $E_X=2B_X+N$, which is clearly an effective class in $X$ such that $E_XN=0$.

Since $E_X$ is an effective divisor such that $E_XB_X<0$ and $E_XN=0$, we conclude that any multiple of $E_X$ is supported on $B_X\cup N$. Thus no positive multiples of $E_X$
can be linearly equivalent to a positive sum of $(-2)$ curves orthogonal to $N$ and hence no multiples of $E_Y$ can be in the effective cone of $Y$.}
\end{rem}

\begin{proposition}\label{prop: Nef cones}
Let us denote by $\Nef(X')$ and $\Nef(Y)$ the nef cones of $X'$ and $Y$
respectively. Then $$\Nef(X')\subseteq \Nef(Y)$$
and the equality holds if and only if there exists no a $(-2)$-curve $B_X$ on $X$ such that $B_XN=1$.
\end{proposition}
\proof We first recall that the lattices $NS(X')$, $NS(Y)$ and $M$
are all isometric.

Let $D_{X'}$ be a nef divisor on $X'$. Then the associated divisor
$D_X$ is a nef divisor on $X$ and $D_XN=0$, by Lemma
\ref{lemma:from X' to X}. This implies, by Lemma \ref{lemma: from
X to Y}, that the corresponding divisor $D_Y$ on $Y$ is nef. So $\Nef(X')\subseteq \Nef(Y)$.

Let $D_Y$ be a nef divisor on $Y$. If there exists no a $(-2)$-curve $B_X$ on $X$ such that $B_XN=1$, then the associated
divisor $D_{X}$ is nef on $X$ and $D_XN=0$ by Lemma
\ref{lemma: from Y to X}. This implies, by Lemma \ref{lemma: from
X to X'} that the corresponding divisor $D_{X'}$ on $X'$ is
nef.

If there exists a $(-2)$-curve $B_X$ on $X$ such that $B_XN=1$, then $\Nef(X')\neq \Nef(Y)$, by Remark \ref{rem: the ccondition on Bx is necessarly}.\endproof

\begin{lemma}\label{lemma: semiample from X to X'}
Let $D_X\in NS(X)$ be a semi-ample divisor (i.e. there exists a positive integer
$m>0$ such that  $|mD_X|$ is without base points) such that
$D_XN=0$. Then the associated divisor $D_{X'}\in NS(X')$ is
semi-ample.
\end{lemma}
\proof Since $D_X$ is semi-ample, there exists an integer $m\in\N$
such that $|mD_X|$ is base points free. In particular $|mD_X|$
does not have fixed component. Since $D_XN=0$, $mD_XN=0$ and thus
$H^0(X,mD_X)$ contains sections which do not pass through any
point of $N$. The images of these sections under the map $\phi$ do
not pass through the point $\phi(N):=P$.

The map $\phi:X\ra X'$ is an isomorphism outside $N$. Let us
consider the sections of $H^0(X,mD_X)$, we call them $s_j$,
$j=1,\ldots, 1+(mD_X)^2/2$. The curves $\phi(s_j)$ are sections of
$H^0(X',mD_{X'})$ and form a basis of this space. Viceversa, let
us consider a section $s'$ in $H^0(X,mD_{X'})$. It is the image of
a curve in $X$ which is in fact the strict transform of $s'$ and
which is clearly a section of $mD_X$. Since
$(mD_X)^2=(mD_{X'})^2$, $h^0(X,mD_X)=h^0(X,mD_{X'})$ and so a
basis of one of these two spaces induces a basis of the other. Let
us assume that $|mD_{X'}|$ has a fixed point $Q$, i.e. all the
sections in $H^0(X', mD_{X'})$ passes through $Q$. We already
observed that $Q\neq P=\phi(N)$, since $P$ is not a base point for
$|mD_{X'}|$. This would imply that $\phi^{-1}(Q)$ is a base point
for $|mD_X|$, which is impossible since $|mD_X|$ is base points
free. So $|mD_{X'}|$ is base points free and thus $D_{X'}$ is
semi-ample.
\endproof
\subsection{Mori Dream Spaces: the surfaces $Y$ and $X'$}
We recall the definition of Mori Dream Space, in the context of
the surfaces (and we recall that any small $\Q$-factorial
modification is an isomorphism if we are considering surfaces).

\begin{definition}{\rm (\cite{HK})} We will call a normal projective surface $X$ a Mori Dream Space provided the following hold:
\begin{enumerate}
\item $X$ is $\Q$-factorial and $\Pic(X)\otimes \Q\equiv
NS(X)\otimes \Q$; \item $\Nef(X)$ is the convex hull of finitely
many semi--ample line bundles.
\end{enumerate}\end{definition}
By \cite[Proposition 2.9]{HK} a variety $X$ such that
$Pic(X)\otimes \Q=NS(X)\otimes \Q$ is a Mori Dream Space if and
only if its Cox ring $\mathcal{R}(X)$ is finitely generated.

\begin{proposition}\label{prop: X' MDS iff Y MDS}
Let us assume that there exists no a $(-2)$-curve $B_X$ on $X$ such that $B_XN=1$. Then the surface $Y$ is a Mori Dream Space if and only if the surface $X'$ is a Mori Dream Space.
\end{proposition}
\proof We recall that every semi-ample divisor is nef. On a K3 surface also the viceversa holds, i.e. a divisor on a K3 surface is semi-ample if and only if it is nef.

We recall that under the assumptions, $\Nef(X')=\Nef(Y)$, by Proposition \ref{prop: Nef cones}.

Let $Y$ be a Mori Dream Space. By definition $\Nef(Y)$ is the
convex hull of a finite set of semi-ample divisors,
$D_{Y}^{(1)},\ldots D_Y^{(r)}$. Then $\Nef(X')$ is the convex hull
of a finite set of nef divisors, which are the divisors
$D_{X'}^{(1)},\ldots D_{X'}^{(r)}$ associated to
$D_{Y}^{(1)},\ldots D_Y^{(r)}$. It remains to show that
$D_{X'}^{(1)},\ldots D_{X'}^{(r)}$ are semi-ample.

The associated divisors on $X$, $D_{X}^{(1)},\ldots D_X^{(r)}$ are
nef and such that $D_X^{i}N=0$ for each $i=1,\ldots r$, by Lemma
\ref{lemma: from Y to X}.  So they are semi-ample divisors on $X$
(which is a K3 surface). By Lemma \ref{lemma: semiample from X to
X'}, the divisors $D_{X'}^{(1)},\ldots D_{X'}^{(r)}$ are
semi-ample on $X'$. Thus $X'$ is a Mori Dream Space.

Let us now assume that $Y$ is not a Mori Dream Space, then $X'$ is
not a Mori Dream Space. Indeed, since $Y$ is a K3 surface, it is
$\Q$-factorial, $Pic(Y)\otimes \Q\equiv NS(Y)\otimes \Q$, and the
nef divisors are semi-ample. So if $Y$ is not a Mori Dream Space,
then $\Nef(Y)$ is not the convex hull of finitely many nef
divisors, which implies that also $\Nef(X')$ is not the convex
hull of finitely many nef divisors, and thus it can not be a Mori
Dream Space.
\endproof

\section{Admissible pairs $(X,X')$}\label{sec: pairs X and X'}

\begin{definition} An admissible pair is a pair of surfaces $(X,X')$ such that $X$ is a K3 surface that is not a Mori Dream Space
and $X'$ is obtained by contracting exactly one extremal ray of
$X$ and is a Mori Dream Space. \end{definition}

Let us consider an admissible pair $(X,X')$. By assumption the
Picard number of $X$ is greater then or equal to 2. Under this
condition, by \cite[Theorem 2]{Kov} the extremal rays of the
effective cone of $X$ are curves with self intersection either 0
or $-2$. Since we are assuming that $X'$ is a surface, the
contraction associated to the extremal ray can not be a fibration
and so we are contracting a $(-2)$-curve.

In order to answer positively to the question \ref{question} it
suffices to find an admissible pair. In this section we describe
how to find the N\'eron--Severi groups of the K3 surfaces $X$ in
admissible pairs (see Condition \ref{condition}) and we give a
list of admissible pairs such that $\rho(X)\geq 10$ (cf.\ Theorem
\ref{theo: the pairs X,X'}).

Let us assume that $X$ and $Y$ are K3 surfaces as in Section
\ref{sec: criterion} (i.e.\ $X$ admits a $(-2)$ curve $N$ such that
the lattice orthogonal to $N$ in $NS(X)$ is the N\'eron--Severi
group of $Y$). We recall that $NS(X')\simeq NS(Y)$, where $X'$ is
obtained by $X$ contracting $N$. By Proposition \ref{prop: X' MDS
iff Y MDS}, if there are no $(-2)$-curves $B_X$ on $X$ such that $B_XN=1$, then the pair of surfaces $(X, X')$ is an admissible pair
if and only if $X$ is not a Mori Dream Space and $Y$ is a Mori
Dream Space. Both $X$ and $Y$ are K3 surfaces and since the K3
surfaces which are Mori Dream Spaces are classified, this provides
a way  to construct pairs $(X,X')$ as required.

Here we summarize some results on K3 surfaces, which will be used
in the following.
\begin{theorem}\label{theorem: when a K3 surface is a Mori Dream Space}{\rm (\cite[Theorems 2.7,2.11, 2.12]{AHL})}
An algebraic K3 surface $S$ is a Mori Dream Space if and only if
the automorphism group $\Aut(S)$ is finite.

In particular, if $\rho(S)=1$, then $S$ is a Mori Dream
Space;\\
if $\rho(S)=2$, then $S$ is a Mori Dream Space if and only if
$NS(S)$ contains at least an element with self intersection either 0 or
$-2$;\\
if $\rho(S)\geq 3$, then $S$ is a  Mori Dream Space if and only if
$NS(S)$ belongs to a finite known list of hyperbolic lattices.
\end{theorem}

The previous Theorem gives a constructive way to produce admissible pairs
$(X,X')$, indeed it suffices to have a K3 surface $X$ with the following properties:
\begin{condition}\label{condition}\begin{itemize}\item $|\Aut(X)|=\infty$; \item there exists a rational curve $N$ on $X$; \item there are no $(-2)$-curves $B_X$ on $X$ such that $B_XN=1$;
\item if $M:=N^{\perp_{NS(X)}}$ and $Y$ is a K3 surface such that
$NS(Y)\simeq M$, then $|\Aut(Y)|<\infty$.\end{itemize}
\end{condition}

\begin{theorem}\label{theo: the pairs X,X'}
Let $(X,X')$ be an admissible  pair of surfaces, then $\rho(X)\geq 3$. 

Let us assume that there exists no a $(-2)$-curve $B_X\subset X$ such that $B_XN=1$. Then $\rho(X)\neq 19$ and if $\rho(X)\geq 4$, then there are finitely many possible choices
for the lattice $NS(X)$.

The complete list of the N\'eron--Severi groups $NS(X)$ of
admissible pairs $(X,X')$ such that there exists no a $(-2)$-curve $B_X\subset X$ with $B_XN=1$ and $\rho(X)\geq 10$ is given in
the table \eqref{table pairs X,X' hight rho}.

If $L$ is as in Table \eqref{table pairs X,X' hight rho} and $X$
is such that $NS(X)\simeq L$,
then there exists a smooth irreducible rational curve $N\subset X$
such that, denoted by $X'$ the surface obtained contracting $N$,
$NS(X')\simeq M$ and thus $(X,X')$ is an admissible pair.

\begin{eqnarray}\label{table pairs X,X' hight rho}\begin{array}{|c|c|c|} \hline
\rho(X)&L\simeq NS(X)&M\simeq NS(X')\\
\hline
20&U\oplus E_8^2\oplus A_1^2&U\oplus

E_8^2\oplus A_1\\
\hline

18&U\oplus E_8\oplus E_7\oplus A_1&U\oplus E_8\oplus E_7\\
\hline

17&U\oplus E_8\oplus D_6\oplus A_1&U\oplus E_8\oplus D_6\\
\hline

16&U\oplus E_8\oplus D_4\oplus A_1^2&U\oplus E_8\oplus D_4\oplus A_1\\
\hline

15&U\oplus E_8\oplus A_1^5&U\oplus D_8\oplus D_4\\
\hline

15&U\oplus E_8\oplus A_1^5&U\oplus E_8\oplus A_1^4\\
\hline

14&U\oplus E_7\oplus A_1^5&U\oplus E_7\oplus A_1^4\\
\hline

14&U\oplus E_8\oplus A_3\oplus A_1&U\oplus E_8\oplus A_3\\
\hline

13&U\oplus D_6\oplus A_1^5&U\oplus D_6\oplus A_1^4\\
\hline

13&U\oplus E_8\oplus A_2\oplus A_1&U\oplus E_8\oplus A_2\\
\hline

12&U\oplus D_4\oplus A_1^6&U\oplus D_4\oplus A_1^5\\
\hline

11&U\oplus E_6\oplus A_2\oplus A_1&U\oplus E_6\oplus A_2\\
\hline

11&U\oplus A_1^9&U\oplus A_1^8\\
\hline

10&U(2)\oplus A_1^8&U(2)\oplus A_1^7\\
\hline

10&U\oplus E_8(2)&U(2)\oplus A_1^7\\
\hline

10&U\oplus A_7\oplus A_1&U\oplus A_7\\
\hline

10&U\oplus D_4\oplus A_3\oplus A_1&U\oplus D_4\oplus A_3\\
\hline

10&U\oplus D_5\oplus A_2\oplus A_1&U\oplus D_5\oplus A_2\\
\hline

10&U\oplus D_7\oplus A_1&U\oplus D_7\\
\hline

10&U\oplus E_6\oplus A_1\oplus A_1&U\oplus E_6\oplus A_1\\
\hline
\end{array}
\end{eqnarray}

\end{theorem}
\proof If $X$ is a K3 surface and it is not a Mori Dream Space
then $\Aut(X)$ is not finite. In particular this implies that
$\rho(X)\geq 2$ and that if  $\rho(X)=2$, then $X$ does not admit
curves with self intersection equal either to 0 or $-2$. Since we
require that $X$ admits a $(-2)$-curve, $\rho(X)\geq 3$.

We now assume that $(X,X')$ is an admissible pair such that there exists no a $(-2)$-curve $B_X$ with $B_XN=1$. Under the latter condition, the fact that $(X,X')$ is an admissible pair is equivalent to the conditions $|\Aut(X)|=\infty$ and $|\Aut(Y)|<\infty$, and hence to a condition on the lattices $L$ and $M$. We first investigate the lattice condition on $L$ and $M$ and after that we discuss of the existence of the curves $N$ and $B_X$.

By hypothesis, $NS(X)$ is an overlattice of finite index of
$M\oplus N$ where $N$ is the class of the $(-2)$-curve contracted
and $M$ is primitively embedded in $NS(X)$. We denote by $r$ the
index of the inclusions $M\oplus N\hookrightarrow NS(X)$. If
$r\neq 1$, then there exists a class $(m+N)/r\in NS(X)$ such that $m\in M$,
$(m+N)/r\not\in M\oplus N$. Clearly this implies that
$N((m+N)/r)=N^2/r=-2/r\in \Z$ and thus $r$ is either 1 or 2. So
$NS(X)$ is either $M\oplus N$ or an overlattice of index 2 of
$M\oplus N$.

The number of hyperbolic lattices $M$ with $\rk(M)\geq 3$ such
that if $Y$ is a K3 surface with $NS(Y)\simeq M$, then
$|\Aut(Y)|<\infty$ is finite. If $(X,X')$ is an admissible pair,
then $NS(Y)\simeq NS(X')\simeq M$, and the admissible choices for
$M$ are finite. The lattice $NS(X)$ is an overlattice of index 1
or 2 of $NS(X')\oplus N\simeq M\oplus A_1$. Since the number
of overlattices of index two of $M\oplus A_1$ is finite (up to
isometries), it follows that the possible choices for $NS(X)$ are
finite if $\rk(M)\geq 3$, i.e. if $\rho(NS(X))\geq 4$.

In order to construct the list of the N\'eron--Severi groups of
the admissible pairs $(X,X')$ such that there exists no a $(-2)$-curve $B_X\subset X$ with $B_XN=1$ (and to exclude the case
$\rho(X)=19$), we check the list of the N\'eron--Severi groups of
the K3 surfaces with a finite group of automorphisms, given in
\cite{K}. We denote by $M$ a lattice in this list. If
$NS(X')\simeq M$, then $X'$ is a Mori Dream Space: Indeed, if  there exists no a $(-2)$-curve $B_X\subset X$ with $B_XN=1$, $X'$ is
a Mori Dream Space if and only if $Y$ is a Mori Dream Space by
Proposition \ref{prop: X' MDS iff Y MDS} and $Y$ is a Mori Dream
Space if and only if $Y$ is a K3 surface with a finite
automorphism group, by Theorem \ref{theorem: when a K3 surface is
a Mori Dream Space}. But $NS(Y)\simeq M$, so $Y$ is a Mori Dream
Space.

So for each $M$ in the list given in \cite{K} we have to construct
the lattice $M\oplus N$ and the overlattices of index 2 of
$M\oplus N$. Each of these lattices is a good candidate to be the
N\'eron--Severi group of $X$.  Let us denote by $L$ one of these
lattices and by $X$ a K3 surface such that $NS(X)\simeq L$. We now
have to check that $X$ is not a Mori Dream Space, so we have to
check that $|\Aut(X)|=\infty$, i.e that $L$ is not in the list
given in \cite{K}. In this way one produces the list of the
possible N\'eron--Severi groups of $X$ and $X'$. Then a geometric
argument can be used in order to show that there exists a model of
$X$ such that the class $N$ represents a smooth irreducible
$(-2)$-curve and so an extremal ray and the analysis of the lattice $L$ allows to conclude that there exists no a $(-2)$-curve $B_X\subset X$ such that $B_XN=1$.

Our first step is to construct the list of the N\`eron--Severi
groups given in Table \eqref{table pairs X,X' hight rho}. We give
all the details for the first lines of the Table, the other cases
are very similar. Let us check the list given in \cite{K} of the
lattices $M$ such that if $Y$ is a K3 surface with $NS(Y)\simeq
M$, then $|\Aut(Y)|<\infty$. If $\rk(M)=19$, then $M\simeq U\oplus
E_8^2\oplus A_1$. So the lattice $L$ is either $M\oplus N\simeq
U\oplus E_8^2\oplus A_1^2$ or an overlattice of index 2 of
$U\oplus E_8^2\oplus A_1^2$. The discriminant group of the lattice
$U\oplus E_8^2\oplus A_1^2$ is $(\Z/2\Z)^2$, so an overlattice of
index two of this lattice is unimodular. But there exits no a
hyperbolic even unimodular lattice of rank 20. So $L\simeq U\oplus
E_8^2\oplus A_1^2$. We now assume $NS(X)\simeq L$. Since $L$ is
not contained in the list given in \cite{K} we conclude that
$|\Aut(X)|=\infty$ and thus $X$ is not a Mori Dream Space. Moreover, since $L\simeq M\oplus \Z N$, all the vectors in $v\in L$ are of type $v:=m+\eta N$, $\eta\in \Z$, $m\in M$, hence $vN\in 2\Z$ and thus there are no divisors $D\in NS(X)$ such that $DN=1$. In particular there is no a $(-2)$-curve $B_X\subset X$ with $B_XN=1$. This
guarantee that if the pair $(X,X')$ is admissible and
$\rho(X)=20$, then $NS(X)\simeq U\oplus E_8^2\oplus A_1^2$ and
$NS(X')\simeq U\oplus E_8^2\oplus A_1$.

Let us consider the lattice $M$ of the list in \cite{K} with
$\rk(M)=18$. In this case $M\simeq U\oplus E_8^2$, which is a
unimodular lattice. A priori $L$ could be either $L\simeq M\oplus
N\simeq U\oplus E_8^2\oplus A_1$ or an overlattice of index 2 of
$M\oplus N\simeq U\oplus E_8^2\oplus A_1$. The discriminant group
of $U\oplus E_8^2\oplus A_1$ is $\Z/2\Z$, so there are no
overlattices of index two of $U\oplus E_8^2\oplus A_1$, and thus
$L\simeq M\oplus N\simeq U\oplus E_8^2\oplus A_1$. Hence, if there
exists an admissible pair $(X,X')$ such that $\rho(X)=19$, then
$NS(X)\simeq U\oplus E_8^2\oplus A_1$. But if $NS(X)\simeq U\oplus
E_8^2\oplus A_1$, then $|\Aut(X)|<\infty$ by \cite{K}, and so $X$
is a Mori Dream Space, by Theorem \ref{theorem: when a K3 surface
is a Mori Dream Space}. We conclude that there are no admissible
pairs $(X,X')$ such that $\rho(X)=19$. The other cases in the Table are similar.

Let $L$ be a lattice given in the second column of Table
\eqref{table pairs X,X' hight rho} and $M$ the corresponding
lattice given in the third column of the same Table. Now we prove
that the generic K3 surface $X$ such that $NS(X)\simeq L$ admits
an extremal contraction (of the curve $N$) such that the corresponding surface $X'$
has the property $NS(X')\simeq M$ and that there is no a $(-2)$-curve $B_X\subset X$ with $B_XN=1$. This implies that $(X,X')$ is
an admissible pair. 

First we observe that for all the pairs $(L,M)$ in Table \ref{table pairs X,X' hight  rho}, except $(U\oplus E_8(-2),U(2)\oplus A_1^7)$,  $L\simeq M\oplus A_1\simeq M\oplus \Z N$ and this implies that for any divisor $D\in NS(X)\simeq  M\oplus \Z$, $DN\in 2\Z$. Thus there exists no a $(-2)$-class $B_X\subset X$ with $B_XN=1$. 
For the pair $(L,M)\simeq (U\oplus E_8(2),U(2)\oplus A_1^7)$, one has to deeply analyze the lattice $L$, which is an overlattice of index 2 of $U(2)\oplus A_1^8$. Denoted by $(u_1,u_2,N_1,\ldots, N_8)$ the basis of $U(2)\oplus A_1^8$, $L$ is obtained adding to this set of divisors the divisor $(\sum_{i=1}^8N_i)/2$. Each $l\in L$ is of the form $l:=a_1u_1+a_2u_2+\sum_{i=1}^8 b_iN_1-k(\sum_{i=1}^8 N_1)/2$, where $a_j,b_i\in \Z$ and $k$ is either 0 or 1. Choosing the curve $N$ to be $N_8$, the condition $lN=lN_8=1$ implies that $k=1$ and $b_8=0$. The condition $l^2=-2$ is now equivalent to $4a_1a_2-2\sum_{i=1}^8b_i^2+2\sum b_i-4=-2$, which is impossible modulo 4. Thus also the pair $(L,M)\simeq (U\oplus E_8(-2),U(2)\oplus A_1^7)$ is such that there exists no a $(-2)$-curve $B_X\subset X$ with $B_XN=1$.

Now it suffices to prove that there exists a
smooth irreducible rational curve which is represented in $NS(X)$
by the class $N$ such that $N^{\perp_{NS(X)}}\simeq M\simeq
NS(X')$. In all the cases but $(L,M)\simeq (U\oplus E_8\oplus A_1^5, U\oplus D_6\oplus D_4), (U(2)\oplus A_1^8, U(2)\oplus A_1^7), (U\oplus E_8(2),U(2)\oplus A_1^7)$, the lattice $L\simeq NS(X)$ is $L\simeq U\oplus R\oplus A_1$ for a certain root lattice $R$ and the lattice $M\simeq NS(X')$ is $M\simeq U\oplus R$.

Since $L\simeq U\oplus R\oplus A_1$, the surface $X$ admits an
elliptic fibration such that the irreducible components of the
reducible fibers which do not intersect the zero section are
represented by the lattice $R\oplus A_1$. In particular, the
lattice $A_1$ which appears as direct summand in the decomposition
$L\simeq U\oplus R\oplus A_1$, represents an irreducible component
of a fiber of type $I_2$ (or $III$). Thus $A_1$ is generated by the
class of a smooth irreducible rational curve. The contraction of
this curve gives the singular surface $X'$, whose N\`eron--Severi
group is naturally identified with $U\oplus R\simeq M$.

In Section \ref{subsect: example with rho=20}, an explicit
equation of the elliptic fibration on $X$ associated to the
decomposition $NS(X)\simeq U\oplus E_8^2\oplus A_1^2$ is provided,
as example.

The existence of the smooth irreducible rational curve $N$ in the
remaining cases (i.e. $(L,M)\simeq (U\oplus E_8\oplus A_1^5, U\oplus D_6\oplus D_4), (U(2)\oplus A_1^8, U(2)\oplus A_1^7), (U\oplus E_8(2),U(2)\oplus A_1^7)$)
is proved in Section \ref{sec: other examples}, (more precisely in
Example \ref{example: same X different X'}, in Section
\ref{subsect: example with rho=10, index 1} and in Section
\ref{subsect: example with rho=10, index 2} respectively).
\endproof

\begin{rem} In Theorem \ref{theo: the pairs X,X'} we proved that
if $(X,X')$ is an admissible pair, there exists no a $(-2)$-curve $B_X\subset X$ such that $B_XN=1$ and $\rho(X)\geq 4$, then the
lattice $NS(X)$ is isometric to a lattice in a finite list of
hyperbolic lattice. On the other hand, if $\rho(X)=3$ then it is
possible to construct infinitely many admissible pairs
$(X_n,X_n')$ such that $n\in\N$ and $NS(X_n)\not\simeq NS(X_{n'})$
if $n\neq n'$. Examples are provided in Proposition \ref{prop:
infinite pairs with K3 with rank 3 not Mori}.
\end{rem}

\begin{rem} There exist examples of both the following cases:
\begin{itemize} \item $(X_1,X_1')$ and $(X_2,X_2')$ are two admissible pairs such that $NS(X_1)\simeq NS(X_2)$, but $NS(X_1')\not\simeq NS(X_2')$
\item $(X_1,X_1')$ and $(X_2,X_2')$ are two admissible pairs such
that $NS(X_1')\simeq NS(X_2')$, but $NS(X_1)\not\simeq NS(X_2)$
\end{itemize} An example of the first case is given in the Table
\ref{table pairs X,X' hight rho}, $\rho(X)=15$ and an example of the second case is given in Table \ref{table pairs X,X' hight rho}, $M\simeq U(2)\oplus A_1^7$. We briefly
describe the geometry of these case in Example \ref{example: same X
different X'}, Section
\ref{subsect: example with rho=10, index 1} and Section
\ref{subsect: example with rho=10, index 2}.

An more exhaustive example of the second case is provided in Section \ref{sec:
geometric example}, where the geometric details are given.
Moreover, an infinite series of examples is presented in
Proposition \ref{prop: infinite pairs with K3 with rank 3 not
Mori}.
\end{rem}

\begin{rem}\label{rem: comditions on M and L to have no the curve BX}
{\rm In the proof of Theorem \ref{theo: the pairs X,X'} we proved and used the following fact: if $L\simeq M\oplus \Z N$, then there exists no a $(-2)$-curve $B_X\subset X$ such that $B_XN=1$. To be more precise, the existence of such a curve implies that $L$ is an overlattice of index 2 of $M\oplus\Z N$ which contains the class $(m+N)/2$, where $m\in M$, $m^2=-6$ and $m/2\in M^{\vee}/M$. Hence sufficient conditions to conclude that there exists no a $(-2)$-curve $B_X\subset X$ with $B_XN=1$ are:\begin{itemize} \item $M\oplus \Z N$ has index 1 in $L$;\item $M$ does not contain vectors with self intersection $-6$;\item $M$ does not contain vectors $m$ of self intersection $-6$ such that $m/2\in M^{\vee}/M$ (for example this is the case if the discriminant quadratic form of $M$ takes value in $\Z$).\end{itemize} 
}
\end{rem}

Since we described the deep relation between the automorphism
group and the property of being a Mori Dream Space, we now
consider the relation between the automorphism groups of the
surfaces involved in our construction. Since $X'$ is obtained from
$X$ by contracting a curves, $\Aut(X')\subset \Aut(X)$. More
precisely, every automorphism $\alpha$ which does not preserve the
curve $N$ does not descend to an automorphism of $X'$ and every
automorphism $\alpha$ which preserves $N$ descend to an
automorphism $\alpha'$ of $X'$. Every automorphism $\alpha'$ of
$X'$ lifts to an automorphism $\alpha$ of $X$ which leaves
invariant the rational curve $N$. Moreover there is also a strong
relation among the automorphism group of $X'$ and $Y$:

\begin{corollary}\label{cor: automorphism group of X' and Y} Let us assume that there is no a $(-2)$-curve $B_X\subset X$ such that $B_XN=1$. Then group of automorphisms of $X'$ is contained in the group of automorphisms of $Y$.\end{corollary}
\proof Let $\alpha'$ be an automorphism of $X'$. It lifts to an
automorphism $\alpha$ of $X$ which preserves $N$. So $\alpha$
induces an effective Hodge isometry of $H^2(X,\Z)$. We denote by
$T_X$ the transcendental lattice of $X$ and we observe that the
isometry induced by $\alpha$ preserves the splitting $M\oplus
N\oplus T_X$. Hence it is a Hodge isometry of $H^2(Y,\Z)$, whose
N\'eron--Severi group is $M$. Moreover, the Nef cone of $Y$ can be
identified with the one of $X'$ by the Proposition \ref{prop: Nef
cones}. Since $\alpha'$ is an automorphism of $X'$, it preserves
the nef cone of $X'$ and so $\alpha^*$ preserves the ample cone of
$Y$. So $\alpha^*$ is an Hodge effective isometry for $Y$ and thus
it is induced by an automorphism $\alpha_Y$ of $Y$.
Let $\alpha_Y$ be an automorphism of $Y$ induced by $\alpha'\in\Aut(X')$ as above. If $\alpha_Y$ is the identity, then $\alpha_Y^*$ is the identity on $H^2(Y,\Z)$ and thus $\alpha^*$ is the identity on $H^2(X,\Z)$. So $\alpha\in\Aut(X)$ is the identity, by Torelli theorem, hence $\alpha'\in \Aut(X')$ is the identity. So the map $\alpha'\mapsto \alpha_Y$ is injective.
\endproof

\section{Two examples}\label{sec: geometric example}

This section is devoted to the geometric description of two
examples. First we consider a lattice $M$ of rank 2 such that if
$Y$ is a K3 surface with $NS(Y)\simeq M$, then $Y$ is a Mori Dream
Space. We describe the geometry and the automorphism group of $Y$
in Section \ref{sec: quartic with a node}: $Y$ admits a model as
quartic hypersurface with a node and a model as double cover of
$\mathbb{P}^2$. For a generic choice of $Y$ the automorphism group
is generated by the cover involution.

In Section \ref{sec: lattices} we construct two different lattices
$L$ of rank 3: one of them is isometric to $M\oplus A_1$, the
other one is the unique even overlattice of index 2 of $M\oplus
A_1$. If the N\'eron--Severi group of a K3 surface is isometric to
one of these two lattices, then the K3 surface is not a Mori Dream
Space, but the contraction of a $(-2)$-curve on it produces a Mori
Dream Space whose N\'eron--Severi group is isometric to $M$, so we
construct the N\`eron--Severi groups of two admissible pairs.

We describe the geometry of the K3 surfaces of these two
admissible pairs in Sections \ref{sec: quartic with two nodes} and
\ref{sec:U+<-4>}: one of them admits a model as quartic in
$\mathbb{P}^3$ with two nodes (which clearly specializes the model
of $Y$) and as double cover of $\mathbb{P}^2$, the other is an
elliptic fibration and admits a model as double cover of the
Hirzebruch surface $\mathbb{F}_4$. We show that the automorphism
group of both these K3 surfaces are infinite, but that several
automorphisms do not descend to automorphism of the contracted
model, which is in fact a Mori Dream Space.

In the following proposition we summarize the results obtained in
this section.
\begin{proposition}
Let $M$ be the lattice $\langle 4\rangle\oplus\langle -2\rangle$.
Let $L$ be an overlattice of $M$ such that:
\begin{itemize}\item There exists a vector $n\in L$, such that $n^2=-2$ and $n^{\perp_L}\simeq M$ \item $L$ admits a primitive embedding in $\Lambda_{K3}$.\end{itemize}
Then $L$ is a hyperbolic even lattice of rank 3 and it is
isometric either to $\langle 4\rangle\oplus\langle -2\rangle\oplus
\langle -2\rangle$ or to $U\oplus \langle -4\rangle$ (see Section
\ref{sec: lattices}). Moreover,
\begin{enumerate}
\item Let $Y$ be a generic K3 surface such that $NS(Y)\simeq M$.
Then $Y$ admits a map $\phi:Y\ra\mathbb{P}^2$ which is a $2:1$
cover branched along a smooth sextic $C_6\subset \mathbb{P}^2$.
There exists a conic $c_2$ which is tangent to $C_6$ in their six
intersection points. The automorphism group of $Y$ is generated by
the cover involution $\iota$, in particular it is finite (see
Section \ref{sec: quartic with a node}). \item Let $X$ be a
generic K3 surface such that $NS(X)\simeq \langle 4\rangle \oplus
\langle -2\rangle\oplus \langle -2\rangle$. Then $X$ admits two
different maps $\phi_i:X\ra(\mathbb{P}^2)_i$, $i=1,2$. Each of
them is a $2:1$ cover branched along a sextic $(C_6)_i\subset
(\mathbb{P}^2)_i$ with one node in the point $(P)_i\in (C_6)_i$.
There exists a conic $(c_2)_i$ such that $(P)_i\not \in (c_2)_i$
and $(c_2)_i$ is tangent to $(C_6)_i$ in their six intersection
points. The automorphism group of $X$ is infinite and contains the
two (non commutative!) cover involutions $\iota_i$.

The map $\phi_1$ (resp. $\phi_2$) contracts exactly one rational
curve, which is $(c_2)_2$ (resp. $(c_2)_1$). Let $X'$ be the
double cover of $\mathbb{P}^2$ branched along $(\mathcal{C}_6)_1$.
Then $((c_2)_2)^{\perp_{NS(X)}}\simeq NS(X')$ is isometric to $M$
and $(X,X')$ is an admissible pair. Moreover $\Aut(X')$ is
generated by the cover involution $\iota_1$ (in particular it is
finite) and the induced effective Hodge isometry coincides with
the one induced by the involution $\iota$ of $Y$ (see Section
\ref{sec: quartic with two nodes}).

\item Let $S$ be a generic K3 surface such that $N(S)\simeq
U\oplus \langle -4\rangle$. Then $S$ admits an elliptic fibration
$\mathcal{E}:S\ra \mathbb{P}^1$ such that $MW(\mathcal{E})=\Z$ is
generated by a section $s_1$ which has a trivial intersection with
the zero section. The automorphism group of $S$ is infinite and
contains the translation by $s_1$ (which is an automorphism of
infinite order).

There exists a $2:1$ map $\varphi:S\ra\mathbb{P}^5$ whose image is
the cone, $\mathcal{C}$, over the twisted rational quartic in
$\mathbb{P}^4$ and let $B\subset\mathcal{C}$ be the branch locus.
The map $\varphi$ contracts exactly one rational curve $s$, which
is a section of the fibration $\mathcal{E}$ and whose image is the
vertex of the cone. Let $S'$ be the double cover of $\mathcal{C}$
branched along $B$. Then $s^{\perp_{NS(S)}}\simeq NS(S')$ is
isometric to $M$ and $(S,S')$ is an admissible pair. Moreover
$\Aut(S')$ is generated by the cover involution (in particular it
is finite) and the induced effective Hodge isometry coincides with
the one induced by the involution $\iota$ of $Y$ (see Section
\ref{sec:U+<-4>}).
\end{enumerate}
\end{proposition}

The following lemma will be essential, since it implies that Proposition \ref{prop: X' MDS iff Y MDS} can be applied for both the pairs $(X, X')$ and $(S, S')$. 
\begin{lemma} The lattice $M\simeq \langle 4\rangle\oplus \langle -2\rangle$ does not contain a vector of length $-6$. Hence $X$ (respectively $S$) does not contain a $(-2)$ curve which intersects $(c_2)_2$ (respectively $s$) in 1 point.\end{lemma}
\proof The quadratic form of $M$ computed on $(x,y)$ is $4x^2-2y^2$. So a vector has length $-6$ if and only if $2x^2-y^2=-3$. This condition implies that $x\equiv 0\mod 3$ and $y\equiv 0\mod 3$. So we write $x=3h$, $y=3k$. The requirement $(3h,3k)$ has length $-6$, is equivalent to $6h^2-3k^2=-1$, which is clearly impossible modulo 3. By Remark \ref{rem: comditions on M and L to have no the curve  BX}, this implies that, denoted by $N$ a vector with self intersection $-2$, neither $M\oplus \Z N$ or an overlattice of index 2 of $M\oplus \Z N$ contains a vector $B$ with length $-2$ such that $BN=1$. In particular this applies to the lattices $NS(X)\simeq M\oplus \Z (c_2)_2$ and $NS(S)$, which is an overlattice of index 2 of $M\oplus \Z s$, and thus $X$ and $S$ are as in the statement.  \endproof

\subsection{Lattice enhancements}\label{sec: lattices}
Let $M$ be the lattice $\langle 4\rangle\oplus \langle-2\rangle$ and let us denote by $\{n_1, n_2\}$ its basis.
The generators of the discriminant group are $M^{\vee}/M\simeq \langle n_1/4, n_2/2\rangle$.

Let $L$ be a lattice such that $L$ admits a primitive embedding in $\Lambda_{K3}$. Since $\Lambda_{K3}$ is an even lattice, $L$ is an even lattice.

We now require that there exists a vector $n\in L$, such that
$n^2=-2$ and $n^{\perp_L}\simeq M$. In particular this implies
that $L$ is an overlattice of finite index $r\in \N$ of $M\oplus
\Z n$. Since $M\oplus \Z n$ is a hyperbolic lattice of rank 3, $L$
is a hyperbolic lattice of rank 3.

We recall that every even hyperbolic lattice of rank less than 11
admits a primitive embedding in $\Lambda_{K3}$. So its enough to
construct all the non isometric even overlattices of $M\oplus \Z
n$ of index $r\in \N$ in order to classify the admissible lattices
$L$.  We already proved that $r$ is either 1 or 2 in proof of
Theorem \ref{theo: the pairs X,X'}.

The first admissible choice for $L$ is $M\oplus \Z n\simeq \langle 4\rangle \oplus \langle-2\rangle\oplus \langle -2\rangle$ (which clearly corresponds to $r=1$).

Now we look for an overlattice of index $r=2$. A $\Q$-basis of
$L$ is given by $n_1$, $n_2$ and $n$. If this is not a $\Z$ basis,
then there exists a vector $w:=(a_1n_1+a_2n_2+a_3n)/2$, $a_i\in
\Z/2\Z$, such that $w\not \in (M\oplus \Z n)$, $wn_1\in\Z$,
$wn_2\in\Z$, $wn\in\Z$ and $ww\in 2\Z$. The condition $ww\in 2\Z$
immediately implies $a_2\equiv a_3\mod 2$. If $a_2$ and $a_3$ are
both even, then $w\equiv n_1/2\mod (M\oplus \Z n)$ and
$(n_1/2)^2=1\not \in 2\Z$. So $a_2\equiv a_3\equiv 1\mod 2$. Again
the condition $ww\in 2\Z$ implies that $a_1\equiv 1\mod 2$. So
there exists a unique overlattice of index $r\neq 1$ of $M\oplus
\Z n$ and  it is the lattice obtained adding to the $\Q$ basis
$\{n_1, n_2,n\}$ the vector $(n_1+n_2+n)/2$. We can now find a
change of bases in order to give a better description of this
overlattice: let us consider the $\Z$-basis
$\{(n_1-n_2-n)/2,(n_1+n_2-n)/2, -n_1+2n\}$. Computing the bilinear
form on this basis we find $U\oplus \langle -4\rangle$, so the
unique even hyperbolic overlattice of index 2 of $M\oplus \Z n$ is
isometric to $U\oplus \langle -4\rangle$.

\subsubsection{Remarks on the orthogonal to $M$}
There exists a unique embedding (up to isometries) of $M$ in $\Lambda_{K3}\simeq U\oplus U\oplus U\oplus E_8\oplus E_8$ which is given by
$$n_1:=\left(\begin{array}{c}1\\2\end{array}\right)\oplus \left(\begin{array}{c}0\\0\end{array}\right)\oplus \left(\begin{array}{c}0\\0\end{array}\right)\oplus \underline 0\oplus \underline 0,\ \ n_2:=\left(\begin{array}{c}0\\0\end{array}\right)\oplus \left(\begin{array}{c}1\\-1\end{array}\right)\oplus \left(\begin{array}{c}0\\0\end{array}\right)\oplus \underline 0\oplus \underline 0.$$
Let us denote by $T$ the lattice $M^{\perp_{\Lambda_{K3}}}$. It is generated by the generators of the third copy of $U$, by the generators of $E_8\oplus E_8$ and by the two vectors:
$$t_1:=\left(\begin{array}{c}-1\\2\end{array}\right)\oplus \left(\begin{array}{c}0\\0\end{array}\right)\oplus \left(\begin{array}{c}0\\0\end{array}\right)\oplus \underline 0\oplus \underline 0,\ \ t_2:=\left(\begin{array}{c}0\\0\end{array}\right)\oplus \left(\begin{array}{c}1\\1\end{array}\right)\oplus \left(\begin{array}{c}0\\0\end{array}\right)\oplus \underline 0\oplus \underline 0.$$
Since $\Lambda_{K3}$ is unimodular, it is an overlattice of index 8 of $M\oplus T$ and indeed $\Lambda_{K3}/(L\oplus M)$ is generated by $(t_1+n_1)/4$ and $(t_2+n_2)/2$.

In order to construct a lattice isometric to $M\oplus \Z n$ we have to identify a vector $n\in T$ with $n^2=-2$. Let us consider the following two possibilities:
\begin{itemize}
\item $n:=\left(\begin{array}{c}0\\0\end{array}\right)\oplus
\left(\begin{array}{c}0\\0\end{array}\right)\oplus
\left(\begin{array}{c}1\\-1\end{array}\right)\oplus \underline
0\oplus \underline 0.$ In this case there is no overlattice of
$L\oplus \Z n$ contained in $\Lambda_{K3}$ and so $L$ is $\langle
4\rangle\oplus \langle -2 \rangle\oplus\langle -2\rangle$. The
orthogonal of such a lattice in $\Lambda_{K3}$  is generated by
the generators of $E_8\oplus E_8$,  by $t_1$, by $t_2$ and by the
vector $\left(\begin{array}{c}0\\0\end{array}\right)\oplus
\left(\begin{array}{c}0\\0\end{array}\right)\oplus
\left(\begin{array}{c}1\\1\end{array}\right)\oplus \underline
0\oplus \underline 0.$ It is isometric to $E_8\oplus E_8\oplus
\langle -4\rangle\oplus \langle2\rangle\oplus \langle 2\rangle.$
\item $n:=t_1+t_2$. In this case the lattice generated by $n_1,
n_2, (n_1+n_2+n)/2$ is primitively embedded in $\Lambda_{K3}$ and
is clearly an overlattice of index 2 of $M\oplus \Z n$. The
orthogonal of such a lattice in $\Lambda_{K3}$  is generated by
the generators of the third copy of $U$, the generators of
$E_8\oplus E_8$ and by $t_1+2t_2$. This lattice is isometric to
$U\oplus E_8\oplus E_8\oplus\langle 4\rangle$.\end{itemize}

\subsection{Some K3 surfaces with Picard number 2}

\subsubsection{Quartic with a node}\label{sec: quartic with a node}

Let us consider a K3 surface $Y$ which admits a model as quartic with exactly one ordinary double point, which is a singularity of type $A_1$. Then there exists a pseudo-ample divisor $H$ in $NS(Y)$ with self intersection 4 and which is orthogonal to a $(-2)$-vector. Indeed  $\varphi_{|H|}:Y\ra \mathbb{P}^3$ contracts exactly one rational curve to the ordinary double point. The class of this curve, called $N_1$, has self intersection $-2$ (since the curve is a rational curve on a K3 surface) and is orthogonal to $H$ (since the curve is contracted by $\varphi_{|H|}$). Hence the lattice $\langle H,N_1\rangle\simeq \langle 4\rangle \oplus\langle -2\rangle$ is primitively embedded in the N\'eron--Severi group of a K3 surface which admits a model as quartic with a singularity of type $A_1$. Moreover the computation of the moduli of the family of the nodal quartics, implies that generically $NS(Y)\simeq \langle 4\rangle\oplus \langle -2\rangle$.

Up to projective transformations, we can assume that the node of
the quartic is in the point
$(1:0:0:0)\in\mathbb{P}^3_{(x_0:x_1:x_2:x_3)}$ and so an equation
for $Y$ is of the form
\begin{equation}\label{eq: quartic with one node}x_0^2F_2(x_1:x_2:x_3)+x_0F_3(x_1:x_2:x_3)+F_4(x_1:x_2:x_3),\end{equation}
where  $F_i$ are generic homogenous polynomials of degree $i$.

Let us consider the projection of the quartic \eqref{eq: quartic
with one node} from $(1:0:0:0)$ to $\mathbb{P}^2_{(x_1:x_2:x_3)}$.
It exhibits $Y$ as double cover of $\mathbb{P}^2_{(x_1:x_2:x_3)}$
branched along the sextic curve
\begin{equation}C_6:=V\left(F_3(x_1:x_2:x_3)^2-4F_2(x_1:x_2:x_3)F_4(x_1:x_2:x_3)\right)\subset\mathbb{P}^2_{(x_1:x_2:x_3)}\end{equation}
The conic $c_2:=V\left(F_2(x_1:x_2:x_3)\right)$ intersects the branch locus in the six points $F_2(x_1:x_2:x_3)=F_3(x_1:x_2:x_3)=0$ and each of them has multiplicity 2.

The projection from the singular point of the quartic is
associated to the divisor $H-N_1$.

The cover involution of the $2:1$ map
$\varphi_{|H-N_1|}:Y\ra\mathbb{P}^2$ is an involution, called
$\iota$, of $Y$ and it does not preserve the symplectic structure
of $Y$ (indeed $Y/\iota$ is birational to $\mathbb{P}^2$). The
class of $H-N_1$ is clearly preserved by the isometry $\iota^*$,
but the class $N_1$ is not. Indeed, since $c_2$ is everywhere
tangent to the branch locus, its inverse image consists of two
disjoint curves. One of them is $N_1$ and the other is $2H-3N_1$.
So $\iota^*$ acts as $-1$ on the transcendental lattice and as
\begin{eqnarray}\label{eq: matrix autom}\left[\begin{array}{rr}3&2\\-4&-3\end{array}\right]\end{eqnarray}
on the basis $\{H, N_1\}$ of the N\'eron--Severi group.

\begin{rem} The K3 surface $Y$ admits two different (equivalent up to automorphism of the surface, but not up to projectivity of $\mathbb{P}^3$) models as singular quartic in $\mathbb{P}^3$. One of them is associated to the divisor $H$, the other one to $\iota^*(H)=3H-4N_1$. %We observe that there are no rational curves on $Y$ which are sent to lines in the model associated to $H$ (since there are no classes in $NS(S)$ which have intersection 1 with $H$). Thus there are no rational curves on $Y$ which are sent to lines in the model associated to $3H-4N_1$ (since there is an automorphism which transforms the model associated to $H$ to the one associated to $3H-4N_1$).
\end{rem}

Since $Y$ has a smooth rational curve and its Picard number is 2,
the automorphism group of $Y$ is finite. To be more precise it is
known that the automorphism group of $Y$ generically coincides with
$\Z/2\Z\simeq \langle \iota\rangle$ (see e.g. \cite{GLP}).

\subsubsection{K3 surfaces with an elliptic fibration}\label{sec: elliptic fibrations}

Let $V$ be a K3 surface and $\mathcal{E}:V\ra\mathbb{P}^1$ be an elliptic fibration (i.e. a fibration in curves of genus 1 which admits at least one section). We will denote by $F$ the class of the fiber of $\mathcal{E}$ and by $s_0$ the class of a given section (called zero section) of $\mathcal{E}$.

Hence $NS(V)\supset \langle F,s_0\rangle$. We obseve that $F^2=0$, $Fs_0=1$ and $s_0^2=-2$ so the bilinear form computed on the basis $\{F,F+s_0\}$ is given by the matrix $U$.

If $V$ is generic among the K3 surfaces admitting an elliptic fibration, then $NS(V)\simeq U$.

On each elliptic curve (so on each smooth fiber of $\mathcal{E}$)
the hyperelliptic involution is well defined. The hyperelliptic
involutions on the fibers glue together giving an involution,
called $h$, of $V$. The quotient by this involution is the
Hirzebruch surface $\mathbb{F}_4$ and the ramification locus
consists of the zero section and of the trisection $t_{10}$ passing
through the 2-torsion of each fiber. The trisection $t_{10}$ is a
curve of genus 10 and is represented by the class $6F+3s_0$.

Since both the class of the section $s_0$ and the class of the
fiber $F$ are preserved by $h$, $h^*$ is the identity on $NS(V)$
and $-1$ on the transcendental lattice.

The divisor $4F+2s_0$ defines a $2:1$ map $\varphi_{|4F+2s_0|}:V\ra \mathcal{C}\subset\mathbb{P}^5$ where $\mathcal{C}$ is a cone over the twisted quartic curve in $\mathbb{P}^4$  (see \cite{SD}). We observe that $\mathcal{C}$ is a model of $\mathbb{F}_4$ obtained contracting the exceptional curve.

The zero section $s_0$ is contracted by $\varphi_{|4F+2s_0|}$ and its image is the vertex of the cone. So the double cover $V\ra\mathcal{C}$ is branched along the image of the trisection, i.e. on $\varphi_{|4F+2s_0|}(t_{10})$. The involution $h$ is exactly the cover involution of this double cover.

Since $V$ contains a rational curve and $\rho(V)=2$, $|\Aut(V)|<\infty$. To be more precise, for a generic choice of $V$, $\Aut(V)=\langle h\rangle$ (see e.g. \cite{GLP}).
\subsection{K3 surfaces with Picard number 3}

\subsubsection{Quartic with two nodes (lattice $L\simeq\langle 4\rangle\oplus \langle-2\rangle\oplus \langle -2\rangle$)}\label{sec: quartic with two nodes}
Let us now consider a K3 surface $X$ such that $NS(X)\simeq \langle 4\rangle \oplus\langle -2\rangle\oplus \langle -2\rangle$. We assume that $X$ is generic among the ones with this property and we denote by $\{H,N_1,N_2\}$ the basis of $NS(X)$.
We can assume that $H$ is a pseudo ample divisor. The map $\varphi_{|H|}:X\ra\mathbb{P}^3$ exhibits $X$ as quartic surface with two double points (the contraction of $N_1$ and $N_2$). Up to projective transformations we can assume that the singular point $\varphi_{|H|}(N_1)$ is $(1:0:0:0)\in\mathbb{P}^3_{(x_0:x_1:x_2:x_3)}$ (as in Section \ref{sec: quartic with a node}) and the singular point $\varphi_{|H|}(N_2)$ is $(0:0:0:1)\in\mathbb{P}^3_{(x_0:x_1:x_2:x_3)}$. So any quartic in this family has the following equation:
\begin{align}\label{eq: quartic with two nodes}\begin{array}{c} x_0^2F_2(x_1:x_2:x_3)+x_0\left(G_3(x_1:x_2)+x_3G_2(x_1:x_2)+x_3^2G_1(x_1:x_2)\right)+\\
+x_3^2H_2(x_1:x_2)+x_3H_3(x_1:x_2)+H_4(x_1:x_2)\end{array}
\end{align}
where $F_i$, $G_i$, $H_i$ are homogenous polynomial of degree $i$.

As in Section \ref{sec: quartic with a node}, we consider the projection from $(1:0:0:0)$ which corresponds to the divisor $H-N_1$. This gives a $2:1$ map $\varphi_{|H-N_1|}:X\ra\mathbb{P}^2_{(x_1:x_2:x_3)}$ which is a double cover branched along the sextic:
$$\begin{array}{c}C_6:=V\left(\left(G_3(x_1:x_2)+x_3G_2(x_1:x_2)+x_3^2G_1(x_1:x_2)\right)^2+\right. \\ \left.-4F_2(x_1:x_2:x_3)\left( x_3^2H_2(x_1:x_2)+x_3H_3(x_1:x_2)+H_4(x_1:x_2)\right)\right).\end{array}$$
The sextic $C_6$ has a unique singular point, which is an ordinary
node, in $(0:0:1)$, i.e. in the point $\varphi_{|H-N_1|}(N_2)$.
Moreover (as in section \ref{sec: quartic with a node}), the conic
$V(F_2(x_1:x_2:x_3))$ intersect $C_6$ in six smooth points and is
tangent to $C_6$ in each of these points. This conic is the image
for $\varphi_{|H-N_1|}$ of the rational curve $N_1$.

The cover involution $\iota_1$ preserves both the class $H-N_1$ (which is the pull back of the hyperplane section of $\mathbb{P}^2_{(x_1:x_2:x_3)}$) and the class $N_2$ (which is sent to the node of the sextic $C_6$ and thus is preserved by the cover involution). Viceversa, the curve corresponding to  $N_1$ is not preserved by $\iota_1$ and $\iota_1^*(N_1)=2H-3N_1$ (see Section \ref{sec: quartic with a node}).
So the involution $\iota_1^*$ on $NS(X)$ is represented, with respect to the basis $\{H,N_1,N_2\}$, by the matrix:
$$\iota_1^*:=\left[\begin{array}{rrr} 3&2&0\\
-4&-3&0\\
0&0&1\end{array}\right].$$

Since the quartic surface \eqref{eq: quartic with two nodes} has two nodes, we can consider two different projections to $\mathbb{P}^2$:  the projection (already considered) from $(1:0:0:0)$ to $\mathbb{P}^2_{(x_1:x_2:x_3)}$, associated to the linear system $|H-N_1|$, which exhibits $X$ as double cover of $\mathbb{P}^2_{(x_1:x_2:x_3)}$ whose cover involution is $\iota_1$; and the projection from $(0:0:0:1)$ to $\mathbb{P}^2_{(x_0:x_1:x_2)}$, which is associated to the linear system $|H-N_2|$.  It exhibits $X$ as double cover of $\mathbb{P}^2_{(x_0:x_1:x_2)}$ and we call the cover involution $\iota_2$. The involution $\iota_2^*$ of $NS(X)$ is represented with respect to the basis $\{H,N_1,N_2\}$ by the matrix
$$\iota_2^*:=\left[\begin{array}{rrr} 3&0&2\\
0&1&0\\
-4&0&-3\end{array}\right].$$

We observe that $\iota_1\iota_2$ is an automorphism of infinite order of $X$, indeed the associated isometry of $NS(X)$ has infinite order and is represented on the basis $\{H,N_1,N_2\}$ by the matrix
$$(\iota_1\iota_2)^*:=\left[\begin{array}{rrr} 9&2&6\\
-12&-3&-8\\
-4&0&-3\end{array}\right].$$

So $\Aut(X)$ is infinite (indeed $\iota_1\iota_2\in\Aut(X)$).

The map $\varphi_{|H-N_1|}$ gives a model of $X$ which contracts
exactly one rational curve, the curve corresponding to the class
$N_2$. Indeed $(H-N_1)^{\perp_{NS(X)}}\simeq \langle -4\rangle
\oplus \langle -2\rangle$ which clearly contains exactly two
vectors with self intersection $-2$ and only one of them is
effective. We observe that $\iota_2^*$ does not preserves the
class $H-N_1$, so it does not descend to an automorphism of
$\varphi_{|H-N_1|}(X)$. Viceversa $\iota_1$ is (by definition) an
automorphism of the model associated to $\varphi_{|H-N_1|}$. The
restriction of $\iota_1^*$ to $N_2^{\perp_{NS(X)}}\simeq NS(Y)$
coincides with the action of $\iota^*$ (where $\iota$ is the
unique non trivial automorphism of $Y$) on $NS(Y)$, given in
\eqref{eq: matrix autom}.

Denoted by $X'$ the (singular) double cover of $\mathbb{P}^2$
branched along $\mathcal{C}_6$, $(X,X')$ is an admissible pair and
$\Aut(X')$ is generated by the cover involution. Indeed, $\iota_1$
induces the cover involution on $X'$, so $\Aut(X')\supset \Z/2\Z$.
By Corollary \ref{cor: automorphism group of X' and Y},
$\Aut(X')\subset \Aut (Y)$ and since $\Aut(Y)=\Z/2\Z$, we deduce
that $\Aut(X')=\Z/2\Z$.

\subsubsection{Elliptic fibrations with non trivial Mordell Weil group (lattices $U\oplus \langle -2d\rangle$)}\label{sec:U+<-2d>}
Here we consider K3 surfaces with Picard number 3 and an elliptic fibration. Since $U$ is primitively embedded in the N\'eron--Severi group of any K3 surface admitting an elliptic fibration and it is a unimodular lattice, the N\'eron--Severi group of a K3 surface with an elliptic fibration and with Picard number 3 is isometric to $U\oplus \langle-2d\rangle$, $d\in \N_{>0}$.

We will denote by $S_d$ a K3 surface such that $NS(S_d)\simeq
U\oplus \langle -2d\rangle$ and we will assume that $S_d$ is
generic among the K3 surfaces with this property. First we discuss
the geometric properties of these surfaces $S_d$ and then we focus
on the case $d=2$.

Let us denote by $\{b_1,b_2,b_3\}$ the basis of $NS(S_d)$ on which the bilinear form is represented by $U\oplus \langle -2d\rangle$. Then we put $F:=b_1$ and $s_0:=b_2-b_1$. The lattice $F^{\perp_{NS(S_d)}}$ consists of the vector $w:=(w_1,w_2,w_3)$ such that $w_2=0$. The bilinear form computed on $w$ is $2w_1w_2-2dw_3^2$, so $F^{\perp_{NS(S_d)}}$ contains a $(-2)$-class if and only if $d=1$.

So if $d=1$, then the elliptic fibration $\varphi_{|F|}:S_d\ra\mathbb{P}^1$ admits exactly one reducible fiber, which is necessarily of type $I_2$, since $S_1$ is generic.

If $d>1$, then the elliptic fibration
$\varphi_{|F|}:S_d\ra\mathbb{P}^1$ has no reducible fibers and
thus, by Shioda--Tate formula, the rank of the Mordell Weil group
is 1. From now on we assume $d>1$ and we denote by $s_1$ a section
of $\varphi_{|F|}:S_d\ra\mathbb{P}^1$ which generates the
Mordell--Weil group. Hence $\{F,s_0,s_1\}$ is a basis of
$NS(S_d)$. It is related to the basis $\{b_1,b_2,b_3\}$ by the
following: $F=b_1$, $s_0=b_2-b_1$, $s_1=(d-1)b_1+b_2+b_3$. It is
immediate to check that the bilinear form computed on this basis
is
$$\left[\begin{array}{rrr}0&1&1\\1&-2&d-2\\1&d-2&-2\end{array}\right].$$
So the K3 surfaces $S_d$, $d>1$, can be geometrically described as
the K3 surfaces admitting an elliptic fibration such that the rank
of the Mordell-Weil group is 1 and the intersection among the zero
section and a generator of the Mordell--Weil group is $d-2$.

We observe that the classes $s_n:=(dn^2-1)b_1+b_2+nb_3$, $n\in \Z$ are the sections of the fibration $\varphi_{|F|}:S_d\ra\mathbb{P}^1$ and there is an isomorphism of groups between the Mordell--Weil group of $\varphi_{|F|}$ and $\Z$ given by $s_n\mapsto n$.
Fixed a value $d$, we have $s_0s_n=dn^2-2$. Moreover $s_is_{i+n}=s_0s_n$, since $s_i$ is the translation (by $s_i$) of $s_0$ and $s_{i+n}$ is the translation (by $s_i$) of $s_n$. It is immediate to check that $s_0s_k>s_0s_h$ if $|k|>|h|$, so $s_0s_1$ is the minimal possible intersection number among two sections of $\varphi_{|F|}:S_d\ra\mathbb{P}^1$.

\subsubsection{The K3 surface $S:=S_2$ ($NS(S)\simeq L\simeq U\oplus \langle -4\rangle$)}\label{sec:U+<-4>}
In case $d=2$, $NS(S_2)\simeq U\oplus\langle -4\rangle$. In the following we will denote by $S$ the surface $S_2$, in order to simplify the notation. The K3 surface $S$ is the generic K3 surface with an elliptic fibration, such that there is a section of infinite order which has a trivial intersection with the zero section.

Our purpose is to describe a map (a geometric model) of $S$ which contracts exactly one rational curve (indeed we already know that there exists a rational curve on $S$ such that the orthogonal to the class of this curve in $NS(S)$ is isometric to the lattice $M\simeq \langle 4\rangle \oplus \langle -2\rangle$, by Section \ref{sec: lattices}).

We consider the map $4F+2s_0$. It exhibits $S$ as double cover of  a cone $\mathcal{C}\subset \mathbb{P}^5$ over the twisted rational curve of degree 4 in $\mathbb{P}^4$, as in Section \ref{sec: elliptic fibrations}. The curve $s_0$ is contracted and it is the unique rational curve which is contracted by this map. Indeed the orthogonal to $4F+2s_0=2b_1+2b_2$ in $NS(S_2)$ is generated by $b_2-b_1=s_0$ and $b_3$. So, if $r$ is the class of a rational curve in $S$ contracted by $4F+2s_0$, then $r=z_1(b_2-b_1)+z_2b_3$, $z_1,z_2\in \Z$; $r^2=-2$; $rF=rb_1\geq 0$. This implies that $z_1=1$ and $z_2=0$, so $r=s_0$.

The family of K3 surfaces which are double covers of
$\mathcal{C}\subset\mathbb{P}^5$ is 18-dimensional (the K3 surface
$V$ described in Section \ref{sec: elliptic fibrations} is a
general member of such a family), but $S$ is the general member of
a 17-dimensional subfamily. Indeed there is a peculiarity in the
branch locus, $\varphi_{|4F+2s_0|}(t_{10})$, of the double cover
$\varphi_{|4F+2s_0|}:S\ra \mathcal{C}$: there exists a curve which
is tangent to $\varphi_{|4F+2s_0|}(t_{10})\subset\mathcal{C}$ in
each of their intersection points. Indeed the class $4F+2s_0$ is
equivalent to $s_1+s_{-1}$ in $NS(S)$. This means that
$\varphi_{|4F+2s_0|}(s_1)=\varphi_{|4F+2s_0|}(s_{-1})$ and thus
the inverse image of $\varphi_{|4F+2s_0|}(s_1)$ splits in the
double cover. So $\varphi_{|4F+2s_0|}(s_1)$ is tangent to
$\varphi_{|4F+2s_0|}(t_{10})$ in all their intersection points,
which are $6=s_1t_{10}$. We observe that since $s_1s_0=0$, the
curve $\varphi_{|4F+2s_0|}(s_1)$ does not pass through the vertex
of $\mathcal{C}$.

The automorphism group of $S$ is surely infinite and contains at
least the following automorphisms: $h$ which is the cover
involution of the double cover $S\ra \mathcal{C}$ (i.e. it is the
hyperelliptic involution of the elliptic fibration
$\varphi_{|F|}:S\ra\mathbb{P}^1$) and the automorphism $T_{s_1}$
which is the translation by the section $s_1$. We observe that:
$h(F)=F$, $h(s_0)=s_0$, $h(s_n)=s_{-n}$ (indeed $h$ switches for
example the sections $s_1$ and $s_{-1}$) and that $T_{s_1}(F)=F$
and $T_{s_1}(s_n)=s_{n+1}$, $n\in \Z$. With respect to the basis
$\{b_1,b_2,b_3\}$ these automorphisms are represented by
$$h=\left[\begin{array}{rrr}1&0&0\\0&1&0\\0&0&-1\end{array}\right]\ \ T_{s_1}=\left[\begin{array}{rrr}1&2&4\\0&1&0\\0&1&1\end{array}\right].$$
It is immediate to check that the class $4F+2s_0=2b_1+2b_2$ is preserved by $h$ but not by $T_{s_1}$, so $T_{s_1}$ does not descend to an automorphism of the singular model $\varphi_{|4F+2s_0|}:S\ra\mathcal{C}$. Indeed the class $s_0$, which is the class of the unique curve contracted by $\varphi_{|4F+2s_0|}$, is not preserved by $T_{s_1}$.

The map $\varphi_{|4F+2s_0|}$ contracts $s_0$ and we observe that
the lattice $s_0^{\perp_{NS(Y_2)}}\simeq M=\langle 4\rangle\oplus
\langle -2\rangle$ is generated by $2b_1+2b_2+b_3$ and
$b_1+b_2+b_3$. The involution $h^*$ restricted to the lattice
$\langle (2b_1+2b_2+b_3), (b_1+b_2+b_3)\rangle$ is represented by
the matrix:
$$\left[\begin{array}{rr}3&2\\-4&-3\end{array}\right],$$
which in fact coincides with the involution $\iota^*$ defined on $M$ by the automorphism $\iota$ of $Y$.

Denoted by $S'$ the (singular) double cover of $\mathcal{C}$
branched along $\varphi_{|4F+2s_0|}(t_{10})$, $(S,S')$ is an
admissible pair and $\Aut(S')$ is generated by the cover
involution. Indeed, $h$ induces the cover involution on $S'$, so
$\Aut(S')\supset \Z/2\Z$. By Corollary \ref{cor: automorphism
group of X' and Y}, $\Aut(S')\subset \Aut (Y)$ and since
$\Aut(Y)=\Z/2\Z$, we deduce that $\Aut(S')=\Z/2\Z$.

\section{Other examples}\label{sec: other examples}
In this section we briefly describe other geometric examples:  We
conclude the proof of the Theorem \ref{theo: the pairs X,X'} in
cases $(L,M)\simeq (U\oplus E_8\oplus A_1^5, U\oplus D_6\oplus D_4), (U(2)\oplus A_1^8, U(2)\oplus A_1^7), (U\oplus E_8(2),U(2)\oplus A_1^7)$, showing the
existence of the smooth irreducible rational curve $N$ which has
to be contracted in order to obtain the Mori Dream Space $X'$. We
also give an explicit equation of an elliptic fibration mentioned
in the proof of the same theorem.

The principal results of this section are geometric descriptions of the following situations:\\
$\bullet$ there exists two admissible pairs $(X_1,X_1')$ and $(X_2,X_2')$ such that $X_1\simeq X_2$ but $X_1'\not\simeq X_2'$ (see Example \ref{example: same X different X'});\\
$\bullet$ there exists two infinite series of admissible pairs
$(S_d, S_d')$ and $(Q_d, Q'_d)$ such that $\rho(Q_d)=\rho(S_d)=3$,
(the minimal possible), and moreover $NS(S_d)\not\simeq NS(Q_d)$
but $NS(S_d')\simeq NS(Q_d')$, see Proposition \ref{prop: infinite
pairs with K3 with rank 3 not Mori}.

\subsection{The K3 surface with $NS(X)\simeq U\oplus E_8^2\oplus
A_1^2$}\label{subsect: example with rho=20} Let $X$ be the K3
surface such that $NS(X)\simeq U\oplus E_8^2\oplus A_1^2$. In the
proof of Theorem \ref{theo: the pairs X,X'} we showed that $X$
admits a smooth irreducible rational curve $N$ which can be
contracted in order to obtain a Mori Dream Space $X'$ and that
this curve is one of the two components of one of the two fibers
of type $I_2$ of a certain elliptic fibration
$\mathcal{E}:X\ra\mathbb{P}^1$. The existence of this fibration
immediately follows from the decomposition of $NS(X)$ in the
direct sum $ U\oplus E_8^2\oplus A_1^2$. Here we give also an
explicit equation of this fibration.

Let us consider the pencil of plane cubics $$V((x^3_1 + x^2_0 x_2
+ x^2_1 x_2) + tx^3_2)\subset \mathbb{P}^2_{(x_0:x_1:x_2)}$$
generated by a triple line $l$ and a smooth cubic $C_3$ which
admits the line $l$ as inflectional tangent. It is well known (and
easy to check) that this pencil induces an elliptic fibration on
the rational surface which is the blow up of $\mathbb{P}^2$ nine
times in the intersection point between $l$ and $C_3$. This
elliptic fibration has one fiber of type $II^*$ (which is the pull
back of $l$) and two fibers of type $I_1$. With the coordinates
$x_2=1$, $y=x_0$ and $x=x_1-1/3$ we immediately obtain the
Weierstrass form
$$y^2=x^3-\frac{1}{3}x+t$$ which has a fiber of type $II^*$ at
infinity and two fibers of type $I_1$ in $t=\pm 2/27$. So the
double cover $\mathbb{P}^1_{(T:S)}\ra\mathbb{P}^1_{(t:s)}$ which
is branched in $t=\pm \frac{2}{27}$, induces the base change
$s:=\frac{27}{4}(S^2-T^2)$ and $t:=\frac{1}{2}(S^2+T^2)$ which
gives the elliptic fibration
$$y^2=x^3-\frac{1}{3}x\left(\frac{27}{4}S^2-\frac{27}{4}T^2\right)^4+\left(\frac{1}{2}S^2+\frac{1}{2}T^2\right)\left(\frac{27}{4}S^2-\frac{27}{4}T^2\right)^5.$$
This is in fact the equation of the unique (up to projective
transformation) elliptic fibration over $\mathbb{P}^1_{(T:S)}$
with reducible fibers $2II^*+2I_2$ and so it is an equation of (a
singular model of) the unique K3 surface $X$ with $NS(X)\simeq
U\oplus E_8^2\oplus A_1^2$. The fibers over $(1:0)$ and $(0:1)$
are of type $I_2$, so each of them consists of two smooth
irreducible rational curves meeting in 2 points. One of these
rational curves meets the zero section of the fibration. The other
can be chosen to be the curve $N$.

\subsection{K3 surfaces with $NS(X)\simeq U(2)\oplus A_1^8$}\label{subsect: example with rho=10, index 1}
We give a geometric description of the K3 surface $X$, which is
generic among the K3 surfaces such that $NS(X)\simeq U(2)\oplus
A_1^8$. We will show that it surely contains a smooth irreducible
rational curve which can be contracted in order to obtain a
singular surface $X'$ whose N\'eron--Severi group is isometric to
$U(2)\oplus A_1^7$. This concludes the proof of Theorem \ref{theo:
the pairs X,X'} in case $NS(X)\simeq U(2)\oplus A_1^8$.

Since there exists a unique even hyperbolic lattice $L$ such that:
$\rk(L)=10$, the discriminant group is $(\Z/2\Z)^{10}$ and the
discriminant form takes values in $\frac{1}{2}\Z$ (and not in
$\Z$), we find $U(2)\oplus A_1^8\simeq \langle 2\rangle\oplus
\langle -2\rangle^9$. So $NS(X)\simeq \langle 2\rangle\oplus
\langle -2\rangle^9$. We can assume that one of the primitive
generators of the sublattice $\langle 2\rangle\hookrightarrow
NS(X)$ is pseudo ample. We denote this divisor by $A$ and we
observe that $\phi_{|A|}:X\ra\mathbb{P}^2$ exhibits $X$ as double
cover of $\mathbb{P}^2$ branched along an irreducible sextic $C_6$
with 9 nodes $P_1,\ldots, P_9$. The double cover of the blow up of
$\mathbb{P}^2$ in these 9 points is a smooth minimal model of $X$
and it contains 9 smooth rational curves which are the double
cover of the 9 exceptional divisors. The classes of these rational
curves are represented by $(-2)$-classes orthogonal to $A$ and
mutually orthogonal. So $X$ admits at least one smooth rational
curve (indeed at least 9) such that the contraction of this curve
produces a singular surface $X'$ whose N\'eron--Severi group is
isometric to $\langle 2\rangle\oplus \langle -2\rangle^8\simeq
U(2)\oplus A_1^7$. A geometric construction of $X'$ is the
following: let us consider the irreducible sextic curve $C_6$ with
9 nodes $P_1,\ldots P_9$ such that $\phi_{|A|}:X\ra\mathbb{P}^2$
is branched along $C_6$. Let us blow up $\mathbb{P}^2$ in the
eight points $P_1,\ldots,P_8$ and let us denote by
$\widetilde{\mathbb{P}^2}$ the surface obtained by this blow up.
Let us denote by $\widetilde{C_6}$ the strict transform of $C_6$
for this blow up. We observe that $\widetilde{C_6}$ has exactly
one singular point. The double cover of $\widetilde{\mathbb{P}^2}$
branched along $\widetilde{C_6}$ is $X'$, indeed the following
diagram commute:
$$\xymatrix{X'\ar[r]^{2:1}&\widetilde{\mathbb{P}^2}\ar[dr]\\
X\ar[r]^{2:1}\ar[u]_{\phi}&\widetilde{\widetilde{\mathbb{P}^2}}\ar[u]\ar[r]&\mathbb{P}^2}
$$
where $\widetilde{\widetilde{\mathbb{P}^2}}$ is the blow up of
$\widetilde{\mathbb{P}}^2$ in the unique singular point of
$\widetilde{C_6}$ and it coincides with the blow up of
$\mathbb{P}^2$ in the nine points $P_1,\ldots, P_9$ and $\phi:X\ra
X'$ is the contraction of the smooth rational curve of $X$ which
is the double cover of the exceptional divisor of the blow up
$\widetilde{\widetilde{\mathbb{P}^2}}\ra\widetilde{\mathbb{P}^2}$.

\subsection{K3 surfaces with $NS(X)\simeq U\oplus E_8(2)$}\label{subsect: example with rho=10, index 2}
We give a geometric description of the K3 surface $X$, which is
generic among the K3 surfaces such that $NS(X)\simeq U\oplus E_8(2)$. We will show that it surely contains a smooth irreducible
rational curve which can be contracted in order to obtain a
singular surface $X'$ whose N\'eron--Severi group is isometric to
$U(2)\oplus A_1^7$. This concludes the proof of Theorem \ref{theo:
the pairs X,X'} in case $NS(X)\simeq  U\oplus E_8(2)$.

Since there exists a unique even hyperbolic lattice $L$ such that:
$\rk(L)=10$, the discriminant group is $(\Z/2\Z)^{8}$ and the
discriminant form takes values in $\Z$, we find that $U\oplus E_8(2)$ is the unique overlattice of index 2 of $U(2)\oplus A_1^8$, and it is generated by the generators of  $U(2)\oplus A_1^8$, $\{u_1,u_2,N_1,\ldots N_8\}$, and by the class $(\sum_{i=1}^8 N_i)/2$.

We can assume that a choice of the primitive
generators of the sublattice $U(2)\hookrightarrow
NS(X)$ consists of nef divisors, and we denote them by $u_1$ and $u_2$. We observe that $\phi_{|u_1+u_2|}:X\ra\mathbb{P}^1\times \mathbb{P}^1\subset\mathbb{P}^3$ exhibits $X$ as double
cover of $\mathbb{P}^1\times\mathbb{P}^1$ branched along a (possibly reducible) curve $\mathcal{C}_{4,4}$ of bidegree $(4,4)$ in $\mathbb{P}^1\times\mathbb{P}^1$ with 8 nodes $P_1=\phi_{|u_1+u_2|}(N_1),\ldots, P_8=\phi_{|u_1+u_2|}(N_8)$. The class of the reduced curve $\phi^{-1}_{u_1+u_2}(\mathcal{C}_{4,4})$ in $NS(X)$ is given by $2u_1+2u_2-\sum_{i=1}^8 N_i$. The curve $\mathcal{C}_{(4,4)}$ is in fact reducible and it is the union of two curves of bidegree $(2,2)$, which corresponds on $X$ to the  class $\left(2u_1+2u_2-\left(\sum_{i=1}^8 N_i\right)\right)/2\in NS(X)$. 
These two components meet exactly in the 8 points $P_i$. The double cover of the blow up of
$\mathbb{P}^1\times \mathbb{P}^1$ in these 8 points is a smooth minimal model of $X$ and it contains 8 smooth rational curves, $N_i$, which are the double cover of the 8 exceptional divisors. The classes of these rational curves are represented by $(-2)$-classes orthogonal to $u_1+u_2$ and mutually orthogonal. So $X$ admits at least one smooth rational
curve (indeed at least 8) such that the contraction of this curve
produces a singular surface $X'$. The N\'eron--Severi group of $X'$ is the orthogonal to an $N_i$, say to $N_8$, in the lattice spanned by  $\{u_1,u_2,N_1,\ldots N_8, (\sum_{i=1}^8 N_i)/2\}$ and it is generated by $\{u_1,u_2,N_1,\ldots N_7\}$. So $NS(X')\simeq
U(2)\oplus A_1^7$. A geometric construction of $X'$ is the
following: let us consider the reducible curve $\mathcal{C}_{4,4}$ with
8 nodes $P_1,\ldots P_8$ such that $\phi_{|u_1+u_2|}:X\ra\mathbb{P}^1\times\mathbb{P}^1\subset \mathbb{P}^3$
is branched along $\mathcal{C}_{4,4}$. Let $\widetilde{\mathbb{P}^1\times\mathbb{P}^1}$ be the blow up $\mathbb{P}^1\times\mathbb{P}^1$ in the
seven points $P_1,\ldots,P_7$.
Let us denote by $\widetilde{\mathcal{C}_{4,4}}$ the strict transform of $\mathcal{C}_{4,4}$
for this blow up. We observe that $\widetilde{\mathcal{C}_{4,4}}$ has exactly one singular point. The double cover of $\widetilde{\mathbb{P}^1\times\mathbb{P}^1}$
branched along $\widetilde{\mathcal{C}_{4,4}}$ is $X'$, indeed the following
diagram commute:
$$\xymatrix{X'\ar[r]^{2:1}&\widetilde{\mathbb{P}^1\times\mathbb{P}^1}\ar[dr]\\
X\ar[r]^{2:1}\ar[u]_{\phi}&\widetilde{\widetilde{\mathbb{P}^1\times\mathbb{P}^1}}\ar[u]\ar[r]&\mathbb{P}^1\times\mathbb{P}^1}
$$
where $\widetilde{\widetilde{\mathbb{P}^1\times\mathbb{P}^1}}$ is the blow up of
$\widetilde{\mathbb{P}^1\times\mathbb{P}^1}^2$ in the unique singular point of
$\widetilde{\mathcal{C}_{4,4}}$ and it coincides with the blow up of
$\mathbb{P}^1\times\mathbb{P}^1$ in the eight points $P_1,\ldots, P_8$ and $\phi:X\ra X'$ is the contraction of the smooth rational curve of $X$ which
is the double cover of the exceptional divisor of the blow up
$\widetilde{\widetilde{\mathbb{P}^1\times\mathbb{P}^1}}\ra\widetilde{\mathbb{P}^1\times\mathbb{P}^1}$.

\subsection{Two admissible pairs $(X,X_1')$ and
$(X,X_2')$ with $X_1'\not\simeq X_2'$ and $\rho(X)=15$} We
consider a K3 surface $X$ whose N\'eron--Severi group is isometric
to $U\oplus D_8\oplus D_4\oplus A_1$, so it is not a Mori Dream
Space (cf. Proposition \ref{theo: the pairs X,X'} and \cite{K}).
We show that it admits two different rational curves $N_1$ and
$N_2$ such that, denoted by $X_i'$ the singular surface
contraction of the curve $N_i$, the two pairs $(X,X_1')$ and
$(X,X_2')$ are both admissible and $X_1'\not\simeq X_2'$. This
gives geometric examples of the cases with $\rho(X)=15$ of the
Table \ref{table pairs X,X' hight rho} and concludes the proof of
Theorem \ref{theo: the pairs X,X'}.
\begin{example}\label{example: same X different X'} Let $X$ be a K3 surface admitting an elliptic fibration $\mathcal{E}:X\ra\mathbb{P}^1$ whose singular fibers
are $I_4^*+I_0^*+I_2+6I_1$ and whose Mordell--Weil group is
trivial. The N\'eron--Severi group of the surface is isomorphic to
$U\oplus D_8\oplus D_4\oplus A_1$ and is generated by the classes
$F,s_0, \Theta_i^{(1)}$, $i=1,\ldots 8$, $\Theta_j^{(2)}$,
$j=1,2,3,4$, $\Theta_1^{(3)}$, where $F$ is the class of the fiber
of the fibration $\mathcal{E}$, $s_0$ is the class of the unique
section of the fibrations, $\Theta_i^{(1)}$ are the eight
components of the fibers of type $I_4^*$ which do not intersect
the zero section, $\Theta_j^{(2)}$  are the 4 components of the
fiber of type $I_0^*$ which do not intersect the zero section,
$\Theta_1^{(3)}$ is the component of the fiber of type $I_2$ which does not intersect the
zero section. The map $\phi:X\ra X_1'$
which contracts the curve $N_1:=\Theta_1^{(3)}$ produces the singular
surface $X_1'$ such that $NS(X_1')\simeq U\oplus D_8\oplus D_4$.

In the following we will assume that the intersection properties
of the components of the fibers of type $I_n^*$ are numbered as
follow:
$$
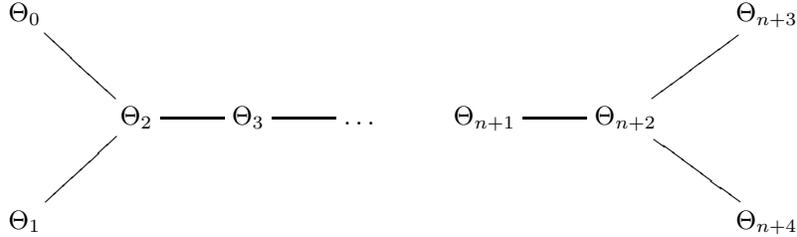
\xymatrix{\Theta_0\ar@{-}[dr]&&&&&&\Theta_{n+3}\\&\Theta_2\ar@{-}[r]&\Theta_3\ar@{-}[r]&\ldots&\Theta_{n+1}\ar@{-}[r]&\Theta_{n+2}\ar@{-}[ur]\ar@{-}[dr]\\
\Theta_1\ar@{-}[ur]&&&&&&\Theta_{n+4}\\}
$$
Let us consider the divisor
$$D:=\Theta_5^{(1)}+2\Theta_4^{(1)}+3\Theta_3^{(1)}+4\Theta_2^{(1)}+5\Theta_0^{(1)}+6s_0+4\Theta_0^{(2)}+2\Theta_2^{(2)}+3\Theta_0^{(3)}$$
on $X$.

We observe that $D^2=0$ and $D$ is an effective divisor. Moreover
the curve $\Theta_6^{(1)}$ is a rational curve such that
$D\Theta_6^{(1)}=1$. Thus, the linear system $|D|$ defines an
elliptic fibration $\phi_{|D|}:X\ra \mathbb{P}^1$. The class of
$D$ is the class of the fiber of the fibration and the divisor $D$
exhibits a reducible fiber of type $II^*$ of this fibration. Let
us denote by $R$ the lattice $\langle
\Theta_5^{(1)},\Theta_4^{(1)},\Theta_3^{(1)},\Theta_2^{(1)},\Theta_0^{(1)},s_0,\Theta_0^{(2)},\Theta_2^{(2)},\Theta_3^{(3)},\Theta_6^{(1)}\rangle$.
Then $R\simeq U\oplus E_8$ and the orthogonal complement of $R$ in
$NS(X)$ is isometric to $A_1^5$ (since the discriminant group of
$NS(X)$ is $(\Z/2\Z)^5$). So the reducible fibers of the elliptic
fibration induced by $|D|$ are $II^*+5I_2$. In particular the
class
$$\Theta_0^{(1)}-\Theta_1^{(1)}+2s_0+2\Theta_0^{(2)}+\Theta_1^{(2)}+2\Theta_2^{(2)}+\Theta_3^{(2)}+\Theta_0^{(3)}$$
is the class of an effective $(-2)$-curve orthogonal to $R$ (and
in fact a bisection of the fibration $\mathcal{E}$) and so it is the
class of one of the components of one of the fibers of type $I_2$ of the fibration $\phi_{|D|}:X\ra \mathbb{P}^1$.
The map which contracts exactly this curve produces a surface
$X_2'$ which is singular in a point and whose N\'eron--Severi group is
$U\oplus E_8\oplus A_1^4$.
\end{example}
\subsection{Infinite admissible pairs $(S_d, S_d')$ with $\rho(S_d)=3$}

Let $S_d$ be a generic K3 surface admitting an elliptic fibration
$\mathcal{E}_d:S_d\ra\mathbb{P}^1$ such that
$MW(\mathcal{E}_d)=\langle s_1\rangle$ and $s_0s_1=d-2$ (described
in Section \ref{sec:U+<-2d>}).

\begin{proposition}\label{prop: infinite pairs with K3 with rank 3 not Mori}  Let $\phi:S_d\ra S'_d$ be the contraction of the curve $s_1$.
Then $NS(S_d')\simeq \langle 2 \rangle\oplus \langle -2d\rangle$. If $d$ is even, then there is no a $(-2)$-curve $B_{S_d}\subset S_d$ with $B_{S_d}s_1=1$.

If $d$ is even and a square, then $(S_d,S_d')$ is an admissible pair. In particular this gives an infinite number of admissible pairs such that the Picard number of the K3 surface is 3 (the minimal possible).

If $d$ is not a square and $d\equiv 0\mod 4$,
then $S'_d$ is not a Mori Dream Space, so $(S_d,S'_d)$ is not an
admissible pair for infinitely many values of
$d$.

Moreover, for almost all the $d$ such that $d$ is a square, the
pair $(Q_d,Q_d')$ is also an admissible if $NS(Q'_d)\simeq
NS(S'_d)$ and $NS(Q_d)\simeq \langle 2 \rangle\oplus \langle
-2d\rangle\oplus \langle -2\rangle$. So we have an infinite number
of admissible pairs $(S_d,S'_d)$ and $(Q_d,Q'_d)$ such that
$NS(S_d)\not\simeq NS(Q_d)$ and $NS(S'_d)\simeq NS(Q'_d)$.
\end{proposition}

\proof We already observed that $S_d$ is not a Mori Dream Space,
since the translation by the section $s_1$ is an automorphism of
infinite order of $S_d$.

Let us assume that there exists a $(-2)$-curve $B_{S_d}\subset S_d$ such that $B_{S_d}s_1=1$. Then there exists a vector $b\in NS(S_d)$ such that $b^2=-2$ and $bs_1=1$. The vector $b$ is of the form $xF+ys_0+zs_1$. So $bs_1=1$ implies $x+(d-2)y-2z=1$, i.e. $x=1+2z-(d-2)y$ and $b^2=-2$ implies $-2y^2-2z^2+2xy+2xz+2(d-2)yz=-2$. These two conditions together give
$y^2+z^2+y+2zy-y^2d+z+1=0$, which is impossible modulo 2 if $d$ is even. We conclude that if $d$ is even there exists no a $(-2)$-curve $B_{S_d}\subset S_d$ such that $B_{S_d}s_1=1$ and thus Proposition \ref{prop: X' MDS iff Y MDS} applies.

The lattice $(s_1)^{\perp_{NS(S_d)}}$ is generated by $\langle
2F+s_1, -dF+s_0-s_1\rangle\simeq \langle 2\rangle\oplus\langle
-2d\rangle $. If $d$ is even, the surface $S_d'$ is a Mori Dream Space if and only if the K3 surface $Y_d$ with N\'eron--Severi group isometric to
$\langle 2\rangle\oplus\langle -2d\rangle$ is a Mori Dream Space.
Since the rank of the lattice is 2, we know that the K3 surface
$Y_d$ is a Mori Dream Space if and only if the lattice $\langle
2\rangle\oplus\langle -2d\rangle$ represents 0 or $-2$.

The quadratic form associated to $\langle 2\rangle\oplus\langle
-2d\rangle$ is $2x^2-2dy^2$.  The form represents the zero if and
only if there exists $(x,y)\in \Z^2$ such that $2x^2-2dy^2=0$ and
represents $-2$ if and only if $2x^2-2dy^2=-2$.

If $d$ is a square, then there exists $b\in\Z$ such that $d=b^2$
and it suffices to chose $(x,y)=(b,1)$. So if $d$ is a square,
then the quadratic form represents zero which implies that $Y_d$
is a Mori Dream Space. So if $d$ is an even square, then $S_d'$ is a Mori Dream Space.

Viceversa if $d$ is not a square, then the quadratic form does not
represent $0$. Let us assume that the quadratic form represents
$-2$. So there exists $(x,y)\in\Z^2$, $x^2-dy^2=-1$. Let us
consider this equation modulo 4 (where the square are either 0 or
1). We have the following possible values for $(x^2,y^2)$ modulo
4, $(0,0)$, $(0,1)$, $(1,0)$, $(1,1)$. The choices $(0,0)$ and
$(1,0)$ give a contradiction. So either $(x^2, y^2)\equiv (0,1)
\mod 4$ and in this case $d\equiv 1\mod 4$ or $(x^2, y^2)\equiv
(1,1) \mod 4$ and in this case $d\equiv 2\mod 4$. Therefore, if
$d\equiv 0 \mod 4$ or $d\equiv 3 \mod 4$, then the quadratic form
does not represent $-2$. If moreover $d$ is not a square, the
quadratic form does not represent also 0.
If $d\equiv 0\mod 4$, then $d$ is even and thus there is no a $(-2)$-curve $B_{S_d}\subset S_d$ such that $B_{S_d}s_1=1$. So if $d\equiv0\mod 4$, then $S_d'$ is a Mori Dream space if and only if $\langle 2\rangle\oplus\langle -2d\rangle$ represents either 0 or $-2$. But if $d\equiv 0\mod 4$ and $d$ is not a square, then $\langle 2\rangle\oplus\langle -2d\rangle$ does not represent neither 0 or $-2$ and thus $S_d'$ in not a Mori Dream Space. 

We observe that $NS(Q_d)\simeq \langle 2\rangle\oplus \langle -2d\rangle\oplus \langle -2\rangle\simeq NS(Y_d)\oplus \langle -2\rangle$ does not contain a vector of length $-2$ which meets the last generator with multiplicity 1, by Remark \ref{rem: comditions on M and L to have no the curve  BX}. So by Proposition \ref{prop: X' MDS iff Y MDS} $Q_d'$ is a Mori Dream Space if and only if $Y_d$ is a Mori Dream Space.

If $d$ is a square, then the K3 surface $Y_d$ is a Mori Dream
Space. As a consequence $Q_d'$ are Mori Dream
Space. Hence the pair $(Q_d,Q_d')$ is an admissible pair if and
only if $Q_d$ is not a Mori Dream Space. Since there are exactly
27 hyperbolic lattices $L$ of rank 3 such that if the
N\'eron--Severi group of a K3 surface is isometric to one of these
lattices, the K3 surface has a finite automorphism group and so is
a Mori Dream Space, we conclude that for almost all the even squares $d$,
$Q_d$ is not a Mori Dream Space and $(Q_d, Q_d')$ is an admissible
pair.

In particular for almost all the $d\in\N$ such that $d$ is an even square, both $(Q_d,Q_d')$ and $(S_d,S_d')$ are admissible pair such that $NS(S_d)\not\simeq NS(Q_d)$ and $NS(S_d')\simeq NS(Q_d')$.
\endproof


\begin{thebibliography}{EKM}
\bibitem[AHL]{AHL} M.\ Artebani, J.\ Hausen, A.\ Laface, {\it On Cox rings of K3-surfaces} Compositio Math. {\bf 146} (2010), 964--998.
\bibitem[GLP]{GLP} F.\ Galluzzi, G.\ Lombardo, C.\ Peters, {\it Automorphs of indefinite binary quadratic forms and K3 surfaces with Picard number 2}, Rendiconti del
Seminario Matematico (Università e Politecnico di Torino), {\bf
68} (2010), 57--77.
\bibitem[L]{L} R.\ Lazarsfeld, Positivity in algebraic geometry. I, Classical setting: line bundles and linear series. {\bf Vol. 48}, Springer-Verlag, Berlin, 2004.
\bibitem[HK]{HK} Y.\ Hu, S.\ Keel {\it Mori Dream Spaces and GIT}, Michigan Math. J. {\bf 48} (2000), 331--348.
\bibitem[Kon]{K} S.\ Kondo, {\it Algebraic K3 surfaces with finite automorphism
group}, Nagoya Math. J. Volume {\bf 116} (1989), 1--15.
\bibitem[Kov]{Kov} S.\ J.\ Kov\'acs , {\it The cone of curves of a K3 surface}. Math. Ann.  {\bf 300}  (1994), 681--691.
\bibitem[O]{O} K.\ Oguiso, {\it A question of Doctor Malte Wandel on automorphisms of the punctural Hilbert schemes of K3 surfaces}, preprint, arXiv:1402.4228
\bibitem[SD]{SD} B.\ Saint-Donat, {\it Projective models of K3
surfaces}, American Journal of Mathematics {\bf 96} (1974),
602--639.
\end{thebibliography}
\end{document}